\documentclass{article}
\pdfoutput1
\usepackage{amssymb,amsmath}
\usepackage{amsfonts}
\usepackage{latexsym}
\newtheorem{Lem}{Lemma}
\newtheorem{Thm}{Theorem}
\newtheorem{Pro}{Proposition}
\newtheorem{Rem}{Remark}
\newtheorem{Cor}{Corollary}
\newtheorem{Exa}{Example}
\newtheorem{Def}{Definition}
\def\C{{\mathbb C}}
\def\R{{\mathbb R}}
\def\N{{\mathbb N}}

\def\H{{\mathbb H}}
 
\def\S{{\mathbb S}}
\def\a{{\mathbf a}}
\def\V{{\mathcal V}}
\newcommand{\Cc}{\mathfrak{C}}
\newcommand{\Kk}{\mathfrak{K}}
\newcommand{\Pp}{\mathfrak{P}}

\def\H{{\mathbb H}}
 
\def\B{{\mathcal B}}
\def\b{{\bf b}}
 
\def\s{{\bf  s}} 
\def\A{{\mathcal{A}}}

\title{Spectrum and Analytic Functional Calculus for  Clifford  Operators via Stem Functions }
\author{Florian-Horia Vasilescu\\
\small Department of Mathematics, University of Lille,\\
\small 59655 Villeneuve d'Ascq, France\\
\small florian.vasilescu@univ-lille.fr}
\date{}

\begin{document}

\maketitle

\begin{abstract} The main purpose of this work is the construction of an analytic functional calculus 
for Clifford operators, which are operators acting on certain  modules over Clifford 
algebras. Unlike in  some preceding works by other authors, we use a spectrum defined in the complex
plane, and certain  stem functions, analytic in neighborhoods of such a spectrum.
The replacement  of  the  slice regular functions, having values in a Clifford algebra,  by analytic 
stem functions becomes possible because of an isomorphism induced by a Cauchy type transform, 
whose existence is proved in the first part of this work.
\end{abstract}
\medskip

{\it Keywords:}  Clifford algebras; stem functions; analytic functional calculus; Clifford and complex spectra.

{\it Mathematics Subject Classification} 2010: 30G35; 30A05; 47A10; 47A60


\section{Introduction}\label{I}

It is unanimously admitted that the analytic functional calculus is a basic tool in the study of linear operators.  While the case of a single operator is settled  by the Riesz-Dunford functional calculus (see for instance \cite{DuSc}), the case of several operators  is more complicate even in the case of commuting ones, because of the non canonical character of  Cauchy type formulas for
analytic functions in several variables, defined in neighborhoods of joint spectra. Nevertheless, a complete construction for  the case of commuting tuples of Banach space operators can be found in the papers \cite{Tay1, Tay2}, and  which was applied, in particular, in investigations related to the local spectral theory (see \cite{Vas1}, \cite{ EsPu},  etc.). 

The case of not necessarily commuting tuples of operators has been approached by several authors, whose
important contributions can be found in works like \cite{Nel}, \cite{Tay3}, and also in \cite{Jef}, \cite{CoSaSt},  where significant
results are obtained by associating the  tuples of operators with some Clifford algebras, which reduces the discussion to the case 
of a single operator. 

In the present paper, we also associate the tupes of operators with some Clifford algebras but, unlike in  \cite{Jef, CoSaSt, CoGaKi},
we replace the  slice regular functions (see Subsection 2.2), defined in some open subsets of Clifford algebras,   by  analytic stem functions (see Definition \ref{stem}), defined in the complex plane. This is a consquence of adopting a different concept of spectrum, that from   \cite{CoSaSt} being 
defined as a  subset of a Clifford algebra,  identified with an Euclidean space, while the spectrum  in this work is defined in the complex plane.  Nevertheless, the 
analytic functional calculi, obtained in two ways, are  shown to be equivalent (see Remark \ref{twofc}). This equivalence is based on the isomorphism of the corresponding spaces of functions which are used in these constructions, and it  is a consequence of some results proved in the first part of this work (see Theorem  \ref{equiv-ons-dom1}).

As in the case of Hamilton's algebra of quaternions, an  important investigation in the  context of  Clifford algebras has been to find a convenient manner to express the  ''analyticity`` of functions, defined on subsets of such algebras. 

A concept of $S$-{\it monogenic function} was introduced in \cite{CoSaSt0}, generalizing that of {\it slice regularity} (see 
\cite{GeSt}) to the framework of Clifford algebres, leading to a significant development sythesized in the first part of \cite{CoSaSt} (which contains a large list of references), whose impact is still actual.   

Unlike in \cite{CoSaSt0}, the basic idea of the present paper to define the regularity of  functions, taking values in a Clifford 
algebra, via an analytic functional calculus, adapting the corresponding results in the Hamilton algebra context  from \cite{Vas2}.

Roughly speaking, and unlike in \cite{GeSt}, an ''$S$-monogenic function`` (which, for simplicity, will be often called in this text ''slice regular``) can and will be obtained by a pointwise application of  the analytic functional calculus with stem functions on a conjugate symmetric open set $U$ in the complex plane, to some elements 
of a given Clifford algebra called {\it  paravectors} (see Subsection 2.1), whose spectra are in $U$,  via the matrix version of  Cauchy's formula (\ref{Cauchy_vect}), with no need of slice derivatives.
In this way, we obtain a whole class of  ''regular functions`` (in fact,  Cauchy transforms of stem functions), eventually shown that
this is  precisely the class of  $S$-monogenic functions, by Theorem  \ref{equiv-ons-dom1}.

Then we initiate some elementary spectral theory for what we call Clifford operators (defined in Subsection 2.3), acting on two-sided modules over a Clifford algebra. A functional calculus with analytic functions is also presented, which happen to be equivalent to the calculus with  
$S$-monogenic functions, as already mentioned above.

The advantage of our approach is its simplicity and   connection with the classical theory, using a spectrum in the complex plane and 
a Riesz-Dunford type of a  kernel. 

The structure of this work is the following. Besides this Introduction, there is a preliminary section, presenting mostly well known concepts and result, necessary in the sequel.  The rest of the work  may virtually be divided  into two parts. A first part (including Sections 3-6)  deals with a description of  an analytic functional calculus for 
stem functions, while the second part (including Sections 7-8)  presents elements of spectral theory for Clifford operators.

This work is inspired by the author's articles \cite{Vas2,Vas3}, dedicated to the quaternionic case, several results having similar proofs to the corresponding one from the 
quoted papers, via  minor modifications. Nevertheless, for the sake of completeness, and for the convenience of the readrs, we give full 
arguments. Some fairly new results are presented in Corollary \ref{intrinsic0}, where the so-called {\it intrinsic functions} (see \cite{CoSaSt},
Definition 3.5.1 of \cite{CoGaKi}, Definition 2.1.2) are described in terms of complex-valued stem functions, and  in  Remark \ref{spmapth}, giving a version of the spectral mapping theorem, which, unlike in  \cite{CoSaSt}, Theorem 3.5.9, appears
 in a classical framework.


\section{Preliminaries} 

 
\subsection{Clifford Algebra and its Complexification}\label{ClAlCo}

First of all, we introduce the concept of (real)  Clifford algebra, in 
a restricted sense (as in \cite{CoSaSt}; see also \cite{HiLo, Jef, Por} etc.). Specifically, in this text, by {\it Clifford algebra}, denoted by ${\Cc}_n$ for a fixed integer $n\ge0$, we mean the unital associative real algebra  having $n+1$ generators, say $e_0=1,e_1,\dots,e_n$, satisfying the relations $e_j^2=-1,\,e_je_k=-e_ke_j$ for all $j,k=1,\ldots,n, j\neq k$. In particular, the real algebra $\R$ is a subalgebra of ${\Cc}_n$. In fact,  ${\Cc}_0=\R,\, {\Cc}_1=\C$, and ${\Cc}_2=\H$, that is, the the real, complex and quaternionic algebras are special cases of Clifford algebras. 
 
Setting $\N_n=\{1,2,\ldots,n\}$, for every ordered subset $J=\{j_1,j_2,\cdots,j_p\}\subset\N_n$, with $j_1<j_2<\cdots<j_p$ and $1\le p\le n$, we put $e_J=e_{j_1}e_{j_2}\cdots e_{j_p}$. 
We use the symbol $J\prec\N_n$ to indicate that $J$ is an oredered set as above. We also assume that $\emptyset\prec\N_n$, 
 $e_\emptyset=1$, $e_{\{j\}}=e_j,\,j=1,\ldots,n$,
and that the family $\{e_J\}_{J\prec\N_n}$ is a basis of the vector space $\Cc_n$. Therefore,  an arbitrary element $\a\in\Cc_n$ can be written as 
\begin{equation}\label{arb_el}
\a=\sum_{J\prec\N_n} a_Je_J,
\end{equation}
where $a_J\in\R$ are uniquely determined for all $J\prec\N_n$. To simplify the notation, we shall put $e_0=e_\emptyset=1$, and $a_0=a_\emptyset$.  The elements of the Clifford algebra $\Cc_n$ will be sometimes called {\it Clifford vectors}, or briefly,
 $Cl$-{\it vectors}.

The linear  subspace of $\Cc_n$ spanned by $\{e_j\}_{j=0}^n$  will be 
denoted by $\Pp_n$. It  plays an important role throughout this work.  The $Cl$-vectors  from  the subspace 
$\Pp_n$, which have the form $\a=a_0+\sum_{k=0}^n a_ke_k$ with 
$a_k\in\R$ for all $k=0,\ldots,n$,  will be called {\it paravectors} (as in \cite{CoSaSt}).  The linear subspace $\Pp_n$  will be often identified with the Euclidean space $\R^{n+1}$, via the linear isomorphism
$$
\Pp_n\ni \sum_{k=0}^n a_ke_k\mapsto(a_0,a_1,\ldots,a_n)\in\R^{n+1}.
$$

 For every 
$\a=\sum_{J\prec\N_n} a_Je_J\in\Cc_n$ we have a decomposition
$\a=\Re(\a)+\Im(\a)$, where $\Re(\a)=a_0$ and $\Im(\a)=
\sum_{\emptyset\neq J\prec\N_n} a_Je_J$, that is, the {\it  real part} and the {\it imaginary part} of the $Cl$-vector  $\a\in\Cc_n$, respectively.

The algebra ${\Cc}_n$ has a norm defined by
\begin{equation}
\vert \a\vert^2=\sum_{J\prec\N_n} a_J^2,
\end{equation}
where $\bf a$ is given by (\ref{arb_el}).

The algebra ${\Cc}_n$ also has an {\it involution} 
${\Cc}_n\ni \a\mapsto \a^*\in {\Cc}_n$, which is defined
via the conditions $e_j^*=-e_j\,(j=1,\ldots,n)$,  $r^*=r\in\R$,
${\bf(ab)^*=b^*a^*}$ for all ${\bf a,b}\in\Cc_n$ (see
\cite{CoSaSt}, Definition 2.1.11).  According to Proposition 2.1.12 from \cite{CoSaSt}, we therefore have $(\a^*)^*=\a,\,{\bf(a+b)^*=a^*+b^*}$. Particularly, if $\a=a_0+\sum_{j=1}^n a_je_j$, then  $\a^*=a_0-\sum_{j=1}^n a_je_j$, for
all $\a\in\Pp_n$.


\begin{Rem}\label{adj}\rm Unlike in \cite{CoSaSt}, we  keep the  name of {\it conjugation} for a different concept  (as in \cite{Vas2}, for the case of quaternions).

To define the {\it conjugation}, we consider the complexification
$\Kk_n=\C\otimes_\R\Cc_n$, identified with the direct sum $\Cc_n+i\Cc_n$. This is a unital algebra with the involution  induced by the involution of $\Cc_n$:
$$
\Kk_n\ni{\bf c=a+{\mathit{i}}b}\mapsto{\bf c^*=a^*-{\mathit{i}}b^*}\in\Kk_n
$$
for all ${\bf a,b}\in\Cc_n$.

The $\R$-linear map
$$
\Kk_n\ni{\bf c=a+{\mathit{i}}b}\mapsto{\bf \bar{c}:=a-{\mathit{i}}b}\in\Kk_n
$$
is a {\it conjugation} on $\Kk_n$, that  is, a real automorphism of 
unital algebras, whose square is the identity.

An important feature 
of this construction is that the elements of the real subalgebra
$\Cc_n$ commute with the complex numbers in the algebra $\Kk_n$. 
Of course, $\bar{\bf a}={ \bf a}$ if and only if ${\bf a}\in\Cc_n$, which is useful
criterion to identify the elements of $\Cc_n$ among those of $\Kk_n$.
\end{Rem}


 \subsection{Slice Regular  Functions}
 
The subspace $\Pp_n$ of paravectors in a Clifford algebra $\Cc_n$
will play an important role in what follows. In fact,  we are particularly interested in functions
defined on open subsets of  $\Pp_n$ (which is  identified with 
$\R^{n+1}$), with values in $\Cc_n$. 

 We shall be dealing with the concept of "slice regularity`` of such functions, which is a form of holomorphy  
(introduced in \cite{CoSaSt0}; see also \cite{CoSaSt}  and several works quoted 
within). If  $\a=a_0+\sum_{k=1}^n a_ke_k\in\Pp_n$ is arbitrary, and so $\a^*=a_0-\sum_{k=1}^n a_ke_k\in\Pp_n$,  we have 
$$
\a\a^*=\a^*\a=\sum_{k=0}^n a_k^2=\vert\a\vert^2.
$$

This shows that the nonnull paravectors $\a\in\Pp_n$ are
invertible in the algebra $\Cc_n$ and, in fact, $\a^{-1}=\vert\a\vert^{-2}\a^*\in\Pp_n$.

For $\Kk_n$-valued functions defined on subsets of $\Pp_n$, the concept of {\it slice regularity} (see also \cite{CoSaSt}) is defined as follows. 

Let $\S_n=\{\mathfrak{s}=\sum_{k=1}^n s_ke_k;\sum_{k=1}^n s_k^2=1\}$, that is, the unit sphere of purely imaginary  
elements of $\Pp_n$. 
It is clear that $\mathfrak{s}^*=-\mathfrak{s}$,  
$\mathfrak{s}^2=-1,\,\mathfrak{s}^{-1}=-\mathfrak{s}$, and
$\vert \mathfrak{s}\vert=1$ for all $\mathfrak{s}\in\S_n$.  Moreover, every nonnull paravector $\a$ can be written as $\a=\Re(\a)+ \vert\a\vert \mathfrak{s}_\a$, with $ \mathfrak{s}_\a=\vert\a\vert^{-1}\Im(\a)\in\S_n$.

Now, let  $\Omega\subset\Pp_n$ be an open set, and let $F:\Omega\mapsto\Kk_n$ be a differentiable function. In the spirit of \cite{CoSaSt0,CoSaSt}, we say that $F$ is {\it right slice regular} on $\Omega$ if for all $\mathfrak{s}\in\S_n$,
$$
 \bar{\partial}_\mathfrak{s}F(x+y\mathfrak{s}):
 =\frac{1}{2}
\left(\frac{\partial}{\partial x}+R_\mathfrak{s}\frac{\partial}{\partial y}\right)F(x+y\mathfrak{s})=0,
$$
on the set $\Omega\cap(\R+\R\mathfrak{s})$,
where $R_\mathfrak{s}$ is the right multiplication of the elements of $\Cc_n$ by $\mathfrak{s}$. 

Unlike in \cite{CoSaSt}, we use the right slice regularity rather than the left one because of a reason to be later explained (see Remark
\ref{left_mult}(1)).
Nevertheless, a left slice regularity can also be defined via the left multiplication of the elements of $\Kk_n$ by elements from $\S_n$. In what follows, the right slice regularity will be simply called {\it slice regularity}. 

As mentioned before,  we are particularly interested by the slice regularity of\break  $\Cc_n$-valued functions, but the concept is valid for $\Kk_n$-valued functions and plays an important role in our discussion.


 \begin{Exa}\label{exa_sr}\rm 
(1) The convergent series of the form $\sum_{m\ge 0}\a_m{\kappa}^m$,  on  balls 
$\{\kappa\in\Pp_n;\vert\kappa\vert<r\}$, with $r>0$ and $\a_m\in\Kk_n$ for all $m\ge0$, are $\Kk_n$-valued slice regular on their domain of definition. Of course, when  $\a_m\in\Cc_n$, such functions are $\Cc_n$-valued  slice regular on their domain of definition.
\end{Exa}


\begin{Rem}\rm The previous discussion is somehow unseemly 
when considering the quaternionic algebra $\H=\Cc_2$. Indeed, in this case
we are mainly interested in $\H$-valued functions, defined on 
subsets of $\H$.  The algebra  $\Cc_2$ is generated by $\{1,e_1,e_2\}$, 
and the vector space $\Pp_2$  generated by this set  is strictly included in 
$\H$. As the algebra $\H$ is also  generated by the set
$\{1,e_1,e_2,e_3\}$, where $e_3=e_1e_2$,  it  is isomorphic to the quotient of the  
algebra  $\Cc_3$ by the two-sided ideal generated by $e_3-e_1e_2$.
Consequently, a separate approach concerning the quaternion algebra
$\H$  (as in \cite{CoSaSt}), rather than an approach in the framework of Clifford algebras, 
 seems to be more appropriate, because it is not a particular case (for $n=2$) of the present approach.
In fact the case $\Cc_2=\H$  is treated in the works \cite{Vas2,Vas3}, strongly related to the present work.

 Note also that  for $\C=\Cc_1$  we have $\Pp_1=\C$.
\end{Rem}


\subsection{Clifford Spaces and Clifford Operators}

Roughly speaking, by a {\it Clifford space} (or a $Cl$-{\it space}) we mean a a two-sided module over a given Clifford algebra $\Cc_n$. A Clifford space is, in particular, a real vector space. A more precise terminology will be given in the following.

Let $\Cc_n$ be a fixed Clifford algebra, and  let  $\V$ be a real vector space. Adapting the framework from \cite{CoSaSt}, the space $\V$ is said to be a {\it right $Cl$-space} if it is a right $\Cc_n$- module, that is, there exists in $\V$   a right multiplication with the elements of $\Cc_n$, such that $x1=x, (x+y)\a=x\a+y\a,\,x(\a+{\bf b})=
x\a+x{\bf b},\, x({\bf a}{\bf b})=(x{\bf a}){\bf b}$ for all $x,y\in\mathcal{V}$ and ${\bf a},{\bf b}\in\Cc_n$.

If $\V$ is a right $Cl$-space which is also a Banach space with the norm $\Vert*\Vert$ such that $\Vert x\a\Vert\le C\Vert x\Vert
\vert\a\vert$ for all $x\in\V$ and $\a\in\Cc_n$, where $C$ is a positive constant, then $\V$ is said to be a {\it right Banch  
Cl-space}. 

In a similar way, one defines the concept of a {\it left $Cl$-space} and that of a {\it left Banach Cl-space}.

 A real (Banach) vector space $\mathcal{V}$ will be said to be a {\it (Banach)  $Cl$-space}  if it is simultaneously a right and a left (Banach) $Cl$-space. 

As for the case of quaternionic operators (see \cite{CoSaSt,Vas3}), it seems to be the framework of  Banch  $Cl$-spaces an appropriate one for the study of some specific  linear  operators, to be defined in the following.

If $\V$ is a real or complex Banach  space, we denote by $\B(\V)$  the algebra of all real or complex bounded  linear operators, respectively. 

Let $\V$ be a fixed Banach $Cl$-space. An operator $T\in\mathcal{B(V)}$ is said  to be {\it right $Cl$-linear} if 
$T(x{\bf a})=T(x){\bf a}$ for all $x\in\mathcal{V}$ and $\a\in\Cc_n$. The set of right $Cl$-linear operators will be 
denoted by $\mathcal{B^{\rm r}(V)}$, which is, in particular, a unital  real Banach algebra.

We shall denote by $R_\a$ (resp. $L_\a$) the right (resp. left)  multiplication operator of the elements of $\V$ with the 
$Cl$-vector $\a\in\Cc_n$.   It is clear that $R_\a,L_\a\in\B(\V)$ for all  $\a\in\Cc_n$. Note also that 
$$
\B^{\rm r}(V)=\{T\in\B(\V);TR_\a=R_\a T,\,\a\in\Cc_n\}.
$$

The elements of the algebra $\B^{\rm r}(\V)$  will be sometimes called {\it  right Clifford (or Cl-) operators}. As we work especially with such operators,
the word ''right`` will be usually omitted.  Note that all operators $L_\a,\, \a\in\Cc_n$, are  $Cl$-operators.

Now, let us consider  the complexification 
$\V_\C$ of $\V$, written as  $\V_\C=\V+i\V$. Because $\V$ is a $\Cc_n$-bimodule, the space $\V_\C$ is actually a two-sided $\Kk_n$-module, via the multiplications 
$$
(\a+i\b)(x+iy)=\a x-\b y+i(\a y+\b x),
(x+iy)(\a+i\b)=x\a-y\b+i(y\a+x\b),
$$ 
for all $\a,\b\in\Cc_n,\,  x,y\in\V$.

For every $T\in \B(\V)$, we consider its natural  ''complex  extension`` to $\V_\C$ given by  $T_\C(x+iy)=Tx+iTy$, for all $x,y\in\V$, 
which is clearly  $\C$-linear,  so $T_\C\in\B(\V_\C)$.  In fact, the map $\B(\V)\ni T\mapsto T_\C\in\B(\V_\C)$ is a unital injective
morphism of real algebras. Moreover, if  $T\in\B^{\rm r}(\V)$, the operator 
$T_\C$ is right $\Kk_n$-linear, that is $T_\C((x+i y)(\a+i\b))=
T_\C(x+iy)(\a+i\b)$ for all $\a+i\b\in\Kk_n,\,x+iy\in\V_\C$, via a direct computation.

 The left and  right
multiplications with  $\a\in\Cc_n$ on $\V_\C$ will be still denoted by $L_\a,R_\a$, respectively, as elements of $\B(\V_\C)$. We set
$$
\mathcal{B}^{\rm r}(\V_\C)=\{S\in \mathcal{B}(\V_\C); SR_\a=R_\a S,\,\a\in\Cc_n\},
$$
which is a unital complex algebra, consisting of all right $\Kk_n$-linear operators on $\V_\C$,  containing all operators $L_\a,
\a\in\Cc_n$. It is easily seen that if $T\in\B^{\rm r}(\V)$, then $T_\C\in \B^{\rm r}(\V_\C)$.


\section{Spectrum of a Paravector}

In the complex algebra $\Kk_n$, we have a natural concept of spectrum,
which can be easily described in the case of paravectors. In fact,
this spectrum is similar to that one introduced in \cite{Vas2} for
quaternions. For this reason, most of the arguments used in \cite{Vas2} apply,
with minor modificatrions, to the actual situation.


\begin{Rem}\label{spectrum}\rm We follow the lines of Remark 1 from \cite{Vas2}.

(1) Because  each paravector commutes in 
$\Kk_n$ with every complex number, we have the identities
\begin{equation}\label{id_spec}
(\lambda-\kappa^*)(\lambda-\kappa)=(\lambda-\kappa)(\lambda-\kappa^*)=\lambda^2-
\lambda(\kappa+\kappa^*)+\vert \kappa\vert^2\in\C,
\end{equation}
for all $\lambda\in\C$ and $\kappa\in\Pp_n$. Therefore, if the complex number 
$\lambda^2-2\lambda\Re(\kappa)+\vert\kappa\vert^2$ is 
nonnull, the element $\lambda-\kappa\in\Kk_n$ is invertible.
Conversely, if the element $\lambda-\kappa\in\Kk_n$ is invertible, assuming 
 $\lambda^2-2\lambda\Re(\kappa)+\vert\kappa\vert^2=0$, the identity 
(\ref{id_spec}) implies $\kappa^*=\lambda$. In other words, $\lambda=\kappa\in\R$, which  is impossible because 
$\lambda-\kappa$ is invertible. Therefore, $\lambda-\kappa$ is invertible if and only if the complex number $\lambda^2-2\lambda\Re(\kappa)+\vert \kappa\vert^2$ is 
nonnull, and therefore 
$$
(\lambda-\kappa)^{-1}=\frac{1}{\lambda^2-2\lambda\Re(\kappa)+\vert \kappa\vert^2}(\lambda-\kappa^*).
$$ 
 Hence, the element      
$\lambda-\kappa\in\Kk_n$ is not  invertible if and only if $\lambda=
\Re(\kappa)\pm i\vert\Im(\kappa)\vert$. In this way, the {\it spectrum} of a paravector $\kappa$ is given by the equality  
$\sigma(\kappa)=\{s_\pm(\kappa)\}$, where 
$s_\pm(\kappa)=\Re(\kappa)\pm i\vert\Im(\kappa)\vert$ are the {\it eigenvalues} of $\kappa$.  

The  argument from above shows, in fact, that if $\lambda-\kappa$
has a left inverse, then $\lambda-\kappa$ is invertible.

(2) As usually, the {\it resolvent set} $\rho(\kappa)$ of a paravector $\kappa\in\Pp_n$ is the set $\C\setminus\sigma(\kappa)$, while the function
 $$\rho(\kappa)\ni\lambda\mapsto(\lambda-\kappa)^{-1}\in\Kk_n$$
 is the {\it resolvent $($function$)$} of $\kappa$, which is a 
$\Kk_n$-valued analytic function on $\rho(\kappa)$. 

(3)  Note that two paravectors $\kappa,
\tau\in\Pp_n$ have the same spectrum if and only if  $\Re(\kappa)=
\Re(\tau)$ and $\vert\Im(\kappa)\vert=\vert\Im(\tau)\vert$.

(4) We recall that $\S_n$ is the unit sphere of purely imaginary paravectors. As already noticed,  every paravector 
$\kappa\in \Pp_n\setminus\R$ can be  written as $\kappa = x + y\mathfrak{s}$, where $x,y$ are real numbers, with $x=\Re(\kappa)$,
$y\in\{\pm \vert \Im(\kappa)\vert\}$, and $\mathfrak{s}\in\{\pm\Im(\kappa)/\vert \Im(\kappa)\vert\}\subset \S_n$. Anyway, we always have 
$\sigma(\kappa)=\{x\pm iy\}$, because $\Im(\kappa)=y\mathfrak{s}$.
Note that, for fixed real numbers $x,y$, the spectrum of $\kappa$ does not depend on $\mathfrak{s}$. Thus, for every $\lambda=u+iv\in\C$ with $u,v\in\R$, we have $\sigma(u+v\mathfrak{s})=\{\lambda,\bar{\lambda}\}$ for all  $\mathfrak{s}\in \S_n$.

(5) The equality $\sigma(\kappa)=\sigma(\tau)$ is clearly an equivalence relation in $\Pp_n$. It is easily seen that 
the equivalence class of an element $\kappa_0=x_0+y_0\mathfrak{s}_0\in\Pp_n$ is given by $\{x_0+y_0\mathfrak{s};\mathfrak{s}\in\mathbb{S}_n\}$. 
 
 (6) Fixing an element $\mathfrak{s}\in\S_n$, 
we have an isometric $\R$-linear map from the complex plane $\C$ into the  space $\Pp_n$, say $\theta_\mathfrak{s}$,  defined by $\theta_\mathfrak{s}(u+iv)=u+v\mathfrak{s},\,u,v\in\R$. For every subset $A\subset\C$, we put
\begin{equation}\label{embedA}
A_\mathfrak{s}=\{x+y\mathfrak{s}; x,y\in\R,x+iy\in A\}=\theta_\mathfrak{s}(A).
\end{equation}
Note that, if $A$ is open in $\C$, then $A_\mathfrak{s}$ is open 
in the $\R$-vector space $\C_\mathfrak{s}\subset\Pp_n$.
 \end{Rem}


The corresponding version of Definition 1 from \cite{Vas2}, adapted to the Cliffordian context, is the following.
 
\begin{Def}\label{NCauchy}\rm  The $\Kk_n$-{\it valued Cauchy kernel} on the open set $\Omega\subset\Pp_n$ is given by
\begin{equation}\label{NvalCk}
\rho(\kappa)\times\Omega\ni(\zeta,\kappa)\mapsto (\zeta-\kappa)^{-1}\in\Kk_n.
\end{equation}
\end{Def}

Recapturing Example 2 from \cite{Vas2}, we get the following.


\begin{Exa}\label{SliceregCk}\rm The $\Kk_n$-{\it valued Cauchy kernel} on the open set $\Omega\subset\Pp_n$ is slice regular. 
Specifically, choosing an arbitrary relatively open set $V\subset\Omega\cap (\R+\R\mathfrak{s})$, and fixing 
$\zeta\in\cap_{\kappa\in V}\rho(\kappa)$, we can write for 
$\kappa\in V$ the equalities
$$
\frac{\partial}{\partial x}(\zeta-x-y\mathfrak{s})^{-1}=
(\zeta-x-y\mathfrak{s})^{-2},
$$
$$
R_\mathfrak{s}\frac{\partial}{\partial y}(\zeta-x-
y\mathfrak{s})^{-1}=-(\zeta-x-y\mathfrak{s})^{-2},
$$
because $\mathfrak{s}^2=-1$, and $\zeta$,  $\mathfrak{s}$ and 
$(\zeta-x-y\mathfrak{s})^{-1}$ commute in $\Kk_n$. Therefore,
$$
\bar{\partial}_\mathfrak{s}((\zeta-\kappa)^{-1})=\bar{\partial}_\mathfrak{s}((\zeta-x-y\mathfrak{s})^{-1})=0,
$$ 
implying the assertion. 
\end{Exa} 

The next comment is inspired by Remark 2 from \cite{Vas2}.

\begin{Rem}\label{eigenval}\rm (1) The discussion about the spectrum
of a paravector can be enlarged, keeping the same framework. 
Specifically, we may regard an element $\kappa\in\Pp_n$ as a left multiplication operator on the algebra $\Kk_n$, denoted by $L_\kappa$, defined by $L_\kappa{\bf a}=\kappa{\bf a}$ for all 
${\bf a}\in\Kk_n$. It is quite clear that $\sigma(L_\kappa)=\sigma(\kappa)$. In this context, we may find the eigenvectors of 
$L_\kappa$, which will be of interest in what follows.  Therefore, we should look for solutions $\nu$  of the equation
$\kappa\nu=s\nu$ in the algebra $\Kk_n$, with $s\in\sigma(\kappa)$. 
Writing $\kappa=\kappa_0+\Im(\kappa)$  with $\kappa_0\in\R$, $s_\pm=\kappa_0\pm i\vert\Im(\kappa)\vert$ and  
$\nu={\bf x}+i{\bf y}$ with ${\bf x},{\bf y}\in\Cc_n$, we obtain 
the equivalent equations
$$
\Im(\kappa){\bf x}=\mp\vert\Im(\kappa)\vert{\bf y},\,\, 
\Im(\kappa){\bf y}=\pm\vert\Im(\kappa)\vert{\bf x},
$$
leading to the solutions
$$
\nu_\pm(\kappa)=\left(1\mp i\frac{\Im(\kappa)}{\vert{\Im(\kappa)}\vert}\right){\bf x}
$$
of the equation $\kappa\nu_\pm=s_\pm\nu_\pm$, where 
${\bf x}\in\Cc_n$ is arbitrary, provided 
${\Im(\kappa)}\neq 0$.  

When $\Im(\kappa)=0$, the  solutions are given by  
$\nu=\kappa_0{\bf a}$, with ${\bf a}\in\Kk_n$ arbitrary.

(2) Every paravector $\mathfrak{s}\in\mathbb{S}_n$ may be associated 
with two elements $\iota_\pm(\mathfrak{s})=(1\mp i\mathfrak{s})/2$ in $\Kk_n$, which are commuting idempotents
such that $\iota_+(\mathfrak{s})+\iota_-(\mathfrak{s})=1$ and
$\iota_+(\mathfrak{s})\iota_-(\mathfrak{s})=0$. For this reason,
setting $\Kk_\pm^{\mathfrak{s}}=\iota_\pm(\mathfrak{s})\Cc_n$, we have a direct sum decomposition 
 $\Kk_n=\Kk_+^{\mathfrak{s}}+\Kk_-^{\mathfrak{s}}$. Explicitly,  if ${\bf a}={\bf u}+i{\bf v}$,
with $\bf u,v\in\Cc_n$, the equation $\iota_+(\mathfrak{s}){\bf x}+
\iota_-(\mathfrak{s}){\bf y}=\bf a$
has the unique solution ${\bf x}={\bf u+\mathfrak{s}v},\,
{\bf y}={\bf u-\mathfrak{s}v}\in\Cc_n$, because  
$\mathfrak{s}^{-1}=-\mathfrak{s}$.

In particular, if $\kappa\in\Pp_n$ and  ${\Im(\kappa)}\neq 0$, setting 
$\mathfrak{s}_{\tilde{\kappa}}=\tilde{\kappa}\vert\tilde{\kappa}\vert
^{-1 }$, where $\tilde{\kappa}=\Im(\kappa)$,
the elements $\iota_\pm(\mathfrak{s}_{\tilde{\kappa}})$ are idempotents, as above. 
\end{Rem}

The next result provides explicit formulas of the spectral 
projections (see \cite{DuSc}, Part I, Section VII.1) associated to the operator $L_{\kappa},\,\kappa\in\Pp_n.$ 
Of course, they are not trivial only if $\kappa\in\Pp_n\setminus\R$
because if $\kappa\in\R$, its spectrum is this real singleton, and the 
only spectral projection is the identity.
\medskip

The statement of the result corresponding to Lemma 1 from \cite{Vas2} sounds like that:


\begin{Lem}\label{spproj} Let $\kappa\in\Pp_n\setminus\R$ be fixed. The spectral projections associated to the eigenvalues $s_\pm(\kappa)$ are given by 
\begin{equation}\label{formproj}
P_\pm(\kappa){\bf a}=\iota_\pm(\mathfrak{s}_{\tilde{\kappa}})\a,\,\,
\a\in\Kk_n.
\end{equation}
Moreover, $P_+(\kappa)P_-(\kappa)=P_-(\kappa)P_+(\kappa)=0$, and 
$P_+(\kappa)+P_-(\kappa)$ is the identity on $\Kk_n$.

When $\kappa\in\R$, the corresponding spectral projection is the identity on $\Kk_n$.
\end{Lem}

{\it Proof.} We shall give a direct argument.  Let us fix a paravector $\kappa$ with $\Im(\kappa)\neq 0$. Next, we  write the general formulas for its spectral projections. Setting $s_\pm= 
s_\pm(\kappa)$, the points $s_+,s_-$ are distinct and  not real. We fix an $r>0$ sufficiently small such that,
setting $D_\pm:=\{\zeta\in \rho(\kappa);\vert\zeta-s_\pm\vert\le r\}$, we have 
$D_\pm \setminus\{s_\pm\}\subset \rho(\kappa)$ and $D_+\cap D_-=\emptyset$. Then 
the spectral projections are given by 
$$
P_\pm(\kappa)=\frac{1}{2\pi i}\int_{\Gamma_\pm} (\zeta-L_\kappa)^{-1}d\zeta
$$
where $\Gamma_\pm$ is the boundary of $D_\pm$.  

Using the equality $L_{\kappa}\nu_\pm(\kappa)=
s_\pm(\kappa)\nu_\pm(\kappa)$
(see Remark \ref{eigenval}), for every $\zeta\in\rho(\kappa)$ and ${\bf h}\in\Cc_n$, we obtain 
$$
(\zeta-L_\kappa)^{-1}(1\mp i\mathfrak{s}_{\tilde{\kappa}}){\bf x}=(\zeta- s_\pm)^{-1}(1\mp  i\mathfrak{s}_{\tilde{\kappa}}){\bf x},
$$
by Remark \ref{eigenval}. Therefore,
$$
P_+(\kappa)(1\mp i\mathfrak{s}_{\tilde{\kappa}}){\bf x}=\frac{1}{2\pi i}\int_{\Gamma_+} (\zeta- s_\pm)^{-1}(1\mp i
\mathfrak{s}_{\tilde{\kappa}}){\bf x} d\zeta, 
$$
and
$$
P_-(\kappa)(1\mp i\mathfrak{s}_{\tilde{\kappa}}){\bf x}=\frac{1}{2\pi i}\int_{\Gamma_-} (\zeta- s_\pm)^{-1}(1\mp i\mathfrak{s}_{\tilde{\kappa}}){\bf x} d\zeta.
$$
Using Cauchy's formula, we deduce that
$$
P_+(\kappa)(1-i\mathfrak{s}_{\tilde{\kappa}}){\bf x}=(1-i\mathfrak{s}_{\tilde{\kappa}}){\bf x},\,P_+(\kappa)(1+i\mathfrak{s}_{\tilde{\kappa}}){\bf x}=0,
 $$
and
$$
P_-(\kappa)(1-i\mathfrak{s}_{\tilde{\kappa}}){\bf x}=0,\,P_-(\kappa)(1+i\mathfrak{s}_{\tilde{\kappa}}){\bf x}=(1+
i\mathfrak{s}_{\tilde{\kappa}}){\bf x},
 $$
for all ${\bf x}\in\Cc_n$.

Fixing an arbitrary element ${\bf a}={\bf u}+i{\bf v}\in\Kk_n$ and 
noticing that $P_\pm(\kappa)$ are $\C$-linear, we clearly obtain
$P_\pm(\kappa){\bf a}=\iota_\pm(\mathfrak{s}_{\tilde{\kappa}})\a,\,\,
\a\in\Kk_n$, which are precisely the formulas (\ref{formproj}) from the statement.

The properties $P_+(\kappa)P_-(\kappa)=P_-(\kappa)P_+(\kappa)=0$, and 
$P_+(\kappa)+P_-(\kappa)$ is the identity on $\Kk_n$ are direct 
consequenceds of the analytic functional calculus associated to a
fixed element $\kappa\in\Pp_n$ in the algebra $\Kk_n$. 
\medskip

By a slight abuse of terminology, the projections $P_\pm(\kappa)$
will be also called the {\it spectral projections} of $\kappa$. 
In fact, as formula (\ref{formproj}) shows, they depend only 
on the imaginary part $\tilde\kappa$ of $\kappa$. 


\begin{Cor}\label{repr1} For every $\kappa\in\Pp_n$ and ${\bf a}\in\Kk_n$ we have
$$
L_\kappa{\bf a}=s_+(\kappa)P_+(\kappa)\a+s_-(\kappa)P_-(\kappa)\a.
$$
\end{Cor}

{\it Proof} Assume first that $\kappa\in\Pp_n\setminus\R$ and ${\bf a}\in\Cc_n$. Then formula (\ref{formproj}) implies that
$$
\a=P_+(\kappa)\a+P_-(\kappa)\a=
\iota_+(\mathfrak{s}_{\tilde{\kappa}})\a+\iota_-(\mathfrak{s}_{\tilde{\kappa}})\a.
$$
Therefore,
$$
L_\kappa\a=\kappa\iota_+(\mathfrak{s}_{\tilde{\kappa}})\a+\kappa\iota_-(\mathfrak{s}_{\tilde{\kappa}})\a=
s_+(\kappa) \iota_+(\mathfrak{s}_{\tilde{\kappa}})\a+ s_-(\kappa)\iota_-(\mathfrak{s}_{\tilde{\kappa}})\a=
$$
$$
s_+(\kappa)P_+(\kappa)\a+s_-(\kappa)P_-(\kappa)\a.
$$
This formula can be extended to elements of $\Kk_n$ because all 
operations within are $\C$-linear.

Of course, this formula also holds for $\kappa=r\in\R$, leading to
$L_r{\bf a}=r{\bf a}$


\section{A General Functional Calculus for\\  Stem Functions}

 In this section, some spaces of $\Kk_n$-valued functions, defined on  subsets of  the complex plane,  will be associated with
spaces of functions, defined on subsets of $\Pp_n$, taking values in the Clifford algebra $\Cc_n$, using spectral methods
and functional calculi. As in \cite{Vas2} (in the context of quaternionic valued functions), this is a general functional calculus, with
arbitrary functions (see Theorem \ref{gfc}).


  \begin{Rem}\label{GFC}\rm As in \cite{Vas2},  Remark 3, an idea from the theory of spectral operators, developed in \cite{DuSc}, Part III, can be applied to define an approprate functional calculus, also useful in the Cliffordian context. Specifically, 
 regarding the algebra  $\Kk_n$ as a (complex) Banach space, and denoting by $\B(\Kk_n)$ the Banach space of all  linear 
operators acting on $\Kk_n$, the operator $L_{\kappa},\,\kappa\in\Pp_n$, (see Remark \ref{eigenval}(1)) is a  particular case of a {\it scalar type} operator, as defined in \cite{DuSc}, Part III, XV.4.1. Its resolution of the identity consists of four projections $\{0,P_\pm(\kappa),$I$\}$, including the null operator $0$ and the identity I, where $P_\pm(\kappa)$ are the spectral projections of
$L_\kappa$, and its integral representation is  given by
$$
L_\kappa=s_+(\kappa)P_+(\kappa)+s_-(\kappa)P_-(\kappa)\in\mathcal{B}(\Kk_n), 
$$
provided by Corollary \ref{repr1}. 

For every function $f:\sigma(\kappa)\mapsto\C$ we may define the operator
$$
f(L_\kappa)=f(s_+(\kappa))P_+(\kappa)+f(s_-(\kappa))P_-({\kappa})\in\mathcal{B}(\Kk_n).
$$
which provides a functional calculus with arbitrary functions on 
the spectrum. More generally, regarding the elements of the algebra
$\Kk_n$ as left multiplication operators on $\Kk_n$, we may extend this formula to functions of the form 
 $F:\sigma({\kappa})\mapsto \B(\Kk_n)$, putting
\begin{equation}\label{GFC1}
F(L_\kappa)=F(s_+(\kappa))P_+(\kappa)+F(s_-(\kappa))P_-(\kappa),
\end{equation}
and keeping this order, which is a ''left functional calculus``, not multiplicative, in general. It is this idea which leads  us to try to define some  $\Cc_n$-valued functions on subsets of $\Pp_n$  associated to certain 
$\Kk_n$-valued functions, defined on some  subsets of $\C$.
\end{Rem}


\begin{Def}\label{consym}\rm (1) A subset $S\subset\C$ is said to be {\it conjugate symmetric} if $\zeta\in S$ if and only if 
$\bar{\zeta}\in S$. 

(2) A subset $A\subset\Pp_n$ is said to be {\it spectrally saturated} 
(as in \cite{Vas2}) if whenever $\sigma(\theta)=\sigma(\kappa)$ for some $\theta\in\Pp_n$ and $\kappa\in A$, we also have $\theta\in A$. 

For an arbitrary $A\subset\Pp_n$, we put $\mathfrak{S}(A)=
\cup_{\kappa\in A}\sigma(\kappa)\subset\C$. 

Conversely, for an arbitrary subset $S\subset\C$, we put
$S_\sigma=\{\kappa\in\Pp_n;\sigma(\kappa)\subset S\}$.
\end{Def}


The discussion from \cite{Vas2}, Remark 4, can be adapted to this context in the following way.

\begin{Rem}\label{consymr}\rm 
(1) If $A\subset\Pp_n$ is spectrally saturated, then $S=\mathfrak{S}(A)$ is conjugate symmetric, and conversely, if $S\subset\C$ is conjugate symmetric, then $S_\sigma$ is spectrally saturated, which
can be easily seen. Moreover, the assignment $S\mapsto S_\sigma$ is injective. Indeed, if $\lambda=u+iv\in S,\,u,v\in\R$, then $\lambda\in\sigma(u+v\mathfrak{s})$ for a fixed $\mathfrak{s}\in\S_n$. If 
$S_\sigma=T_\sigma$ for some $T\subset\C$, we must have $\kappa=u+v\mathfrak{s}\in T_\sigma$. Therefore $\sigma(u+v\mathfrak{s})
\subset T$, implying $\lambda\in T$, and so $S\subset T$. Clearly,
we also have $T\subset S$. 

Similarly, the assignment $A\mapsto\mathfrak{S}(A)$ is injective and  $A=S_\sigma$ if and only if $S=\mathfrak{S}(A)$. These two assertions are left to the reader.

(2)   If $\Omega\subset\Pp_n$ is an open spectraly saturated set, then $\mathfrak{S}(\Omega)\subset\C$ is open. To see that, let $\lambda_0=u_0+iv_0\in\mathfrak{S}(\Omega)$ be fixed, with $u_0,v_0\in\R$,  and let 
$\kappa_0=u_0+v_0\mathfrak{s}$, where 
$\mathfrak{s}\in\S_n$ is also fixed. Because $\Omega$ is spectrally 
saturated, we must have $\kappa_0\in\Omega$. As the set $\Omega\cap\C_\mathfrak{s}$ is relatively open, there is a positive number
$r$ such that the open set 
$$
\{\kappa=u+v\mathfrak{s}; u,v\in\R,\vert\kappa-\kappa_0\vert<r\}
$$
is in $\Omega\cap\C_\mathfrak{s}$. Therefore, the
set of the points $\lambda=u+iv$, satisfying $\vert\lambda-\lambda_0\vert<r$ is in $\mathfrak{S}(\Omega)$, implying that it is open. 

Conversely,
if $U\subset\C$ is open and conjugate symmetric,
the set $U_\sigma$ is also open via the upper semi-continuity of the spectrum (see \cite{DuSc}, Part I, Lemma VII.6.3.).

An important particular case is when $U={\mathbb D}_r:=\{\zeta\in\C;\vert\zeta\vert<r\}$, for some $r>0$. Then\
$U_\sigma=\{\kappa\in\Pp_n;\vert\kappa\vert<r\}$. Indeed, if $\vert\kappa\vert
< r$ and  $\theta$ has the property $\sigma(\kappa)=\sigma(\theta)$, from the equality $\{\Re(\kappa)\pm i\vert\Im(\kappa)\vert\}=
\{\Re(\theta)\pm i\vert\Im(\theta)\vert\}$ it follows that 
$\vert\theta\vert<r$.

 (3) A subset $\Omega\subset\Cc_n$ is said to 
be {\it axially symmetric} if for every $\kappa_0=u_0+v_0\mathfrak{s}_0\in\Omega$ with $u_0,v_0\in\R$ and $\mathfrak{s}_0\in\S_n$, we also
have $\kappa=u_0+v_0\mathfrak{s}\in\Omega$ for all $\mathfrak{s}\in\S_n$. This concept is introduced in \cite{CoSaSt}, Definition 2.2.17. In fact, we have the following.


\begin{Lem}\label{UH_tildeU} A subset $\Omega\subset\Cc_n$ is axially
 symmetric if and only if it is spectrally saturated. 
\end{Lem}

The assertion follows easily from the fact that the equality $\sigma(\kappa)=\sigma(\tau)$ is an equivalence relation in
$\Pp_n$  (see Remark \ref{spectrum}(5)).

Nevertheless, we continue to use the expression  ''spectrally saturated set`` to designate an ''axially symmetric set``, because the
former name is more compatible with our spectral approach.  
\end{Rem}

As noticed in Remark \ref{adj}, the algebra $\Kk_n$ is endowed with a conjugation given
by $\bar{\bf a}={\bf b}-i{\bf c}$, when ${\bf a}={\bf b}+i{\bf c}$, with ${\bf b},{\bf c}\in\Cc_n$. Note also that, because $\C$
is a subalgebra of $\Kk_n$, the conjugation of $\Kk_n$ restricted to
$\C$ is precisely the usual complex conjugation.
\medskip 

The next definition has an old origin, seemingly going back to \cite{Fue}.


\begin{Def}\label{stem}\rm Let $U\subset\C$ be conjugate symmetric, and let $F:U\mapsto\Kk_n$. We say that $F$ is a ($\Kk_n$-{\it valued}) {\it stem function} if  $F(\bar{\lambda})=\overline{F(\lambda)}$ for all 
$\lambda\in U$. 
\end{Def}
\medskip

For an arbitrary conjugate symmetric subset $U\subset\C$, we put 
\begin{equation}
\mathcal{S}(U,\Kk_n)=\{F:U\mapsto\Kk_n;F(\bar{\zeta})=\overline{F(\zeta)},\zeta\in U\},
\end{equation}
that is, the $\R$-vector space of all $\Kk_n$-valued stem functions on $U$.  
Replacing $\Kk_n$ by $\C$, we denote by $\mathcal{S}(U)$ 
the real algebra of all $\C$-valued stem functions, which is an 
$\R$-subalgebra in $\mathcal{S}(U,\Kk_n)$. In addition, the space
$\mathcal{S}(U,\Kk_n)$ is a two-sided $\mathcal{S}(U)$-module.  
 

The following definitions adapts, to our context, Definition 4 from \cite{Vas2}.

\begin{Def}\label{funccalc} Let  $U\subset\C$ be  conjugate symmetric.
For every $F:U\mapsto\Kk_n$ and all $\kappa\in U_\sigma$ we define 
a function $F_\sigma:U_\sigma\mapsto\Kk_n$, via the assignment 
\begin{equation}\label{GFC2}
U_\sigma\setminus\R\ni\kappa\mapsto F_\sigma(\kappa)=F(s_+(\kappa))\iota_+(\mathfrak{s}_{\tilde\kappa})+F(s_-(\kappa))\iota_-(\mathfrak{s}_{\tilde\kappa}) \in\Kk_n,
\end{equation} 
where $\tilde\kappa=\Im(\kappa),\,\mathfrak{s}_{\tilde\kappa}=\vert\tilde\kappa\vert^{-1 }\tilde\kappa$, and  
$\iota_\pm(\mathfrak{s}_{\tilde\kappa})=2^{-1}(1\mp i\mathfrak{s}_{\tilde\kappa})$, and $F_\sigma(r)=F(r)$, if $r\in U_\sigma\cap\R$.
\end{Def} 

Formula (\ref{GFC2}) is strongly related to formula  
(\ref{GFC1}) because the spectral projections $P_\pm(\kappa)$ are the left  multiplications defined by $2^{-1}(1\mp i\mathfrak{s}_{\tilde\kappa})$ respectively, via formula (\ref{formproj}).

The next result is a version of Theorem 2 from \cite{Vas2}. 


\begin{Thm}\label{sym_spec} Let $U\subset\C$ be a conjugate symmetric subset, and let $F:U\mapsto\Kk_n$.
The element $F_\sigma(\kappa)$ belongs to $\Cc_n$ for all $\kappa\in U_\sigma$ if and only if  $F\in\mathcal{S}(U,\Kk_n)$.
\end{Thm}

{\it Proof.}\,  We first assume that $F_\sigma(\kappa)\in\Cc_n$ for all $\kappa\in U_\sigma$. We fix a point $\zeta\in U$, supposing that $\Im\zeta>0$. Then we choose a paravector 
$\kappa\in U_\sigma$ with $\sigma(\kappa)=\{\zeta,\bar{\zeta}\}$. 
Therefore, $s_+(\kappa)=\zeta$ and $s_-(\kappa)=\bar{\zeta}$.
We write $F(\zeta)=F_{+1}+iF_{+2}, F(\bar{\zeta})=F_{-1}+iF_{-2}$, 
with $F_{+1},F_{+2}, F_{-1},F_{-2}\in\Cc_n$. 
According to (\ref{GFC2}), we must have 
$$
2F(\kappa)=F_{+1}+F_{+2}\mathfrak{s}_{\tilde\kappa}+F_{-1}-F_{-2}\mathfrak{s}_{\tilde\kappa}+i(-F_{+1}\mathfrak{s}_{\tilde\kappa}
+F_{+2}+F_{-2}+F_{-1}\mathfrak{s}_{\tilde\kappa}),
$$
so
$$
2\overline{F(\kappa)}=F_{+1}+F_{+2}\mathfrak{s}_{\tilde\kappa}+F_{-1}-F_{-2}\mathfrak{s}_{\tilde\kappa}+i(F_{+1}\mathfrak{s}_{\tilde\kappa}-F_{+2}-F_{-2}-F_{-1}\mathfrak{s}_{\tilde\kappa}).
$$
Because $F(\kappa)=\overline{F(\kappa)}$, we must have
$$
-F_{+1}\mathfrak{s}_{\tilde\kappa}
+F_{+2}+F_{-2}+F_{-1}\mathfrak{s}_{\tilde\kappa}=F_{+1}\mathfrak{s}_{\tilde\kappa}-F_{+2}-F_{-2}-F_{-1}\mathfrak{s}_{\tilde\kappa},
$$
which is equivalent to    
$$
F_{+2}+F_{-2}=(F_{+1}-F_{-1})\mathfrak{s}_{\tilde\kappa}
$$  
  
Choosing another paravector $\theta\in\Pp_n$ with $s_+(\theta)=\zeta$
but with $\mathfrak{s}_{\tilde\kappa}\neq\mathfrak{s}_{\tilde\theta}$, 
we obtain
$$
(F_{+1}-F_{-1})(\mathfrak{s}_{\tilde\kappa}-\mathfrak{s}_{\tilde\theta})=0.
$$
Because the paravector $\mathfrak{s}_{\tilde\kappa}-\mathfrak{s}_{\tilde\theta}$ is nonnull, it must be invertible, and so 
$F_{+1}=F_{-1}$, implying $F_{+2}=-F_{-2}$, meaning that 
$\overline{F(\zeta)}=F(\bar{\zeta})$. 

If $\Im\zeta=0$, so $\zeta=x_0\in\R$, taking $\kappa=x_0$, we have $\sigma(\kappa)=\{x_0\}$, and $F(\kappa)=F(x_0)$ is a 
real number.

If $\Im\zeta<0$, applying the above argument to $\bar{\zeta}$ we 
obtain $\overline{F(\bar{\zeta})}=F({\zeta})$. Consequently, $F$
is a stem function on $U$. 

Conversely, if  $\overline{F(\zeta)}=F(\bar{\zeta})$ for all $\zeta\in U$, choosing a $\kappa\in\Pp_n$ with $\Im\kappa\neq0$,
and setting $\zeta=s_+(\kappa)$, we obtain from (\ref{GFC2})
the equality
$$
2F_\sigma(\kappa)=F(\zeta)+F(\bar{\zeta})-i(F(\zeta)-F(\bar{\zeta})\mathfrak{s}_{\tilde\kappa}.
$$
Therefore,
$$
2\overline{F_\sigma(\kappa)}=F(\bar{\zeta})+F(\zeta)+i(F(\bar{\zeta})-F(\zeta))\mathfrak{s}_{\tilde\kappa},
$$
showing that  $F_\sigma(\kappa)\in\Cc_n$ for all $\kappa\in U_\sigma$,
because the case $\kappa=r\in\R$ is obvious.  


\begin{Cor}\label{intrinsic0} Let $U\subset\C$ be a conjugate symmetric subset, and let $f:U\mapsto\C$. The following conditions are equivalent;

$(1)$  $f\in{\mathcal S}(U)$;

$(2)$  $f_\sigma(\kappa)$  belongs to $\C_\mathfrak{s}$, and $f_\sigma(\kappa^*)=f_\sigma(\kappa)^*$  for all $\kappa\in U_\sigma\cap\C_\mathfrak{s}$, where $\C_\mathfrak{s}=\{u+v\mathfrak{s};u,v\in\R\}$, and $\mathfrak{s}\in\S_n$.
\end{Cor}

{\it Proof.} That $f\in{\mathcal S}(U)$ implies  $f_\sigma(\kappa)$  belongs to $\Cc_n$ for all $\kappa\in U_\sigma$
is a direct consequence of  Theorem \ref{sym_spec}. More precisely, in this case actually (2) holds true.
 To get (2), let us first choose a function  $f\in{\mathcal S}(U)$. It follows from Definition
\ref{funccalc} that
$$
f_\sigma(\kappa)=f(s_+(\kappa))\iota_+(\mathfrak{s}_{\tilde\kappa})+
\overline{f(s_+(\kappa))}\iota_-(\mathfrak{s}_{\tilde\kappa})=
$$
$$
\Re(f(s_+(\kappa))+\Im(f(s_+(\kappa))\mathfrak{s}_{\tilde\kappa}\in\C_{\mathfrak{s}_{\tilde\kappa}}
$$
where  $\kappa\in U_\sigma,\,\tilde\kappa=\Im(\kappa),\,\mathfrak{s}_{\tilde\kappa}=\tilde\kappa\vert\tilde\kappa\vert^{-1 }$,
 and  $\iota_\pm(\mathfrak{s}_{\tilde\kappa})=2^{-1}(1\mp i\mathfrak{s}_{\tilde\kappa})$. Because $s_+(\kappa^*)=s_+(\kappa)$, and $\tilde{\kappa}^*=-\tilde{\kappa}$, we clearly have  $f_\sigma(\kappa^*)=f_\sigma(\kappa)^*$  for all $\kappa\in U_\sigma$. 

Conversely, let $g:U_\sigma\mapsto\Pp_n$ be such that  $g(\kappa)\in\C_\mathfrak{s}$, and $g(\kappa^*)=g(\kappa)^*$  for all $\kappa\in U_\sigma\cap\C_\mathfrak{s}$. We shall look
for a function  $f\in{\mathcal S}(U)$ such that $f_\sigma=g$. Fixing the  points $z_\pm=x\pm iy\in U$  with $x,y(\neq0)\in 
\R$, and a paravector $\mathfrak{s}\in\S_n$, we set ${\kappa}_\pm=x\pm y \mathfrak{s}$, so $\kappa_\pm^*=\kappa_\mp$,
 $s_+( {\kappa}_\pm)=x+i\vert y\vert$, and $s_-( {\kappa}_\pm)=x-i\vert y\vert$. Then we define
$$
2f(z_+)=g({\kappa}_+)(1-i\mathfrak{s})+g({\kappa}_-)(1+i\mathfrak{s}),
$$
$$
2f(z_-)=g(\kappa_+)(1+i\mathfrak{s})+g({\kappa}_-)(1-i\mathfrak{s}).
$$
Because $g(\kappa_+)\in\C_\mathfrak{s}$, there are $u,v\in\R$ such that $g(\kappa_+)=u+v\mathfrak{s}$, and thus, $g(\kappa_-)=u-v\mathfrak{s}$. Therefore, $f(z_+)=u+iv\in\C$. Similarly,  $f(z_-)=u-iv\in\C$. This shows, in fact, that 
$f:U\mapsto\C$, and $f(z)=\overline{f(\bar{z})}$. Using the equations from above, we derive easily that
$$
2g({\kappa}_+)=f(z_+)(1-i\mathfrak{s})+f(z_-)(1+i\mathfrak{s}),
$$
$$
2g({\kappa}_-)=f(z_+)((1+i\mathfrak{s})+(z_-)(1-i\mathfrak{s}).
$$
implying the equality
 $g(\kappa)=f_\sigma(\kappa)$ for all $\kappa\in U_\sigma$, which concludes the proof.


\begin{Rem}\label{zeros}\rm  This is a description of the zeros of the functions obtained via Theorem \ref{sym_spec},
corresponding to Remark 5 from \cite{Vas2}.

 Let $U\subset\C$ be a conjugate symmetric set and let $F\in\mathcal{S}(U,\Kk_n)$ be arbitrary. We can easily describe the zeros of $F_\sigma$. 
Indeed, if $F_\sigma(\kappa)=F(s_+(\kappa))\iota_+({\tilde\kappa})+F(s_-(\kappa))\iota_-({\tilde\kappa})=0$, we must have 
$F(s_+(\kappa))\iota_+({\tilde\kappa})=0$ and 
$F(s_-(\kappa))\iota_-({\tilde\kappa})=0$, via a direct manipulation with the idempotents $\iota_\pm({\tilde\kappa})$. In other words, we must have $F(s_\pm(\kappa))=\pm iF(s_\pm(\kappa))\mathfrak{s}_{\tilde\kappa}$. As in the previous proof, choosing another paravector
${\theta}$ with $s_+(\kappa)\in\sigma(\theta)$ and
$\tilde{\kappa}\neq\tilde{{\theta}}$, we obtain
$F(s_+(\kappa))(\mathfrak{s}_{\tilde{\kappa}}-\mathfrak{s}_{\tilde{\theta}})=0$. Therefore, $F(s_+(\kappa))=0$ because 
$ \mathfrak{s}_{\tilde{\kappa}}-\mathfrak{s}_{\tilde{\theta}}$ is
invertible. Similarly, $F(s_-(\kappa))=0$. Conqequently, setting $\mathcal{Z}(F):=\{\lambda\in U;F(\lambda)=0\}$,
and $\mathcal{Z}(F_\sigma):=\{\kappa\in U_\sigma;F_\sigma(\kappa)=0\}$,
we must have
$$
\mathcal{Z}(F_\sigma)=\{\kappa\in U_\sigma;\sigma(\kappa)\subset\mathcal{Z}(F)\}. 
$$
\end{Rem}

For every subset  $\Omega\subset\Pp_n$, we denote by $\mathcal{F}(\Omega,\Cc_n)$ the set of all $\Cc_n$-valued functions on $\Omega$. Let also
\begin{equation}\label{intrinsic}
\mathcal{IF}(\Omega,\Cc_n)=\{g:\mathcal{F}(\Omega,\Cc_n); g(\kappa^*)=g(\kappa)^*\in\C_\mathfrak{s}, \kappa\in\Omega\cap\C_\mathfrak{s},\, \mathfrak{s}\in\S_n\},
\end{equation}
which is a unital commutative subalgebra of the algebra $\mathcal{F}(\Omega,\Cc_n)$. The functions from the space 
$\mathcal{IF}(\Omega,\Pp_n)$ are similar to those called {\it intrinsic functions}, appering in \cite{CoSaSt}, Definition 3.5.1, or in \cite{CoGaKi},
Definition 2.1.2.

The next result provides  a  {\it $\Cc_n$-valued general functional calculus} for arbitrary stem functions (for the quaternionic case
see \cite{Vas2}, Theorem 2). 


\begin{Thm}\label{gfc} Let $\Omega\subset\Pp_n$ be a spectrally
saturated set, and let $U=\mathfrak{S}(\Omega)$. The map
$$
{\mathcal S}(U,\Kk_n)\ni F\mapsto F_\sigma\in\mathcal{F}(\Omega,\Cc_n)
$$ 
is $\R$-linear, injective, and has the property $(Ff)_\sigma=F_\sigma f_\sigma$ for all
$F\in{\mathcal S}(U,\Kk_n)$ and $f\in{\mathcal S}(U)$. Moreover, the
restricted map 
$$
{\mathcal S}(U)\ni f\mapsto f_\sigma\in\mathcal{IF}(\Omega,\Cc_n)  
$$ 
is unital and multiplicative.
\end{Thm}

{\it Proof.}\, The map $F\mapsto F_\sigma$ is clearly $\R$-linear.
The injectivity of this map follows from Remark \ref{zeros}. 
Note also that
$$
F_\sigma(\kappa) f_\sigma(\kappa)=(F(s_+(\kappa))\iota_+(\mathfrak{s}_{\tilde\kappa})+F(s_-(\kappa))\iota_-(\mathfrak{s}_{\tilde\kappa}))\times 
$$
$$
(f(s_+(\kappa))\iota_+(\mathfrak{s}_{\tilde\kappa})+f(s_-(\kappa))\iota_-(\mathfrak{s}_{\tilde\kappa})= 
$$
$$
(Ff)(s_+(\kappa))\iota_+(\mathfrak{s}_{\tilde\kappa})+(Ff)(s_-(\kappa))\iota_-(\mathfrak{s}_{\tilde\kappa})=(Ff)_\sigma (\kappa),
$$ 
because $f$ is complex valued, and using the properties of the idempotents $\iota_\pm(\mathfrak{s}_{\tilde\kappa})$
In particular, this computation shows that if 
$f,g\in{\mathcal S}(U)$, and so $f_\sigma,g_\sigma\in\mathcal{IF}(\Omega,\Cc_n)$  by Corollary \ref{intrinsic0},  we have  $(fg)_\sigma=f_\sigma g_\sigma= 
g_\sigma f_\sigma $, thus the map $f\mapsto f_\sigma$ is multiplicative. It is also clearly unital. 
 

\section{A Cauchy Transform in the Clifford\\  Algebra  Context}

Having  the $\Kk_n$-valued Cauchy kernel (see Definition \ref{NCauchy}), we may introduce a concept of a Cauchy transform  (as in \cite{Vas2} in the quaternionic context),  whose some useful  properties will be exibited  in this section.

The frequent use of versions of the Cauchy formula is simplified by adopting the following definition. Let $U\subset\C$ be open.  An open subset 
$\Delta\subset U$ will be called a {\it Cauchy domain} (in $U$) if 
$\Delta\subset\bar{\Delta}\subset U$ and the boundary $\partial\Delta$ of $\Delta$ consists of a finite family of closed curves, piecewise smooth, positively oriented. A  Cauchy domain is bounded but not necessarily connected.
\medskip

For a given open set $U\subset\C$, we denote by $\mathcal{O}(U,\Kk_n)$ the complex algebra  of all $\Kk_n$-valued analytic functions on $U$. 

If $U\subset\C$ is also  conjugate symmetric, let 
$\mathcal{O}_s(U,\Kk_n)$ be the real subalgebra of $\mathcal{O}(U,\Kk_n)$ consisting of all stem functions from $\mathcal{O}(U,\Kk_n)$. 

Because $\C\subset\Kk_n$, we have $\mathcal{O}(U)\subset\mathcal{O}(U,\Kk_n)$, where $\mathcal{O}(U)$ is the complex algebra  of all 
complex-valued analytic functions on the open set $U$. Similarly, when $U\subset\C$ is also conjugate symmetric, $\mathcal{O}_s(U)\subset\mathcal{O}_s(U,\Kk_n)$, where $\mathcal{O}_s(U)$ is the 
real subalgebra consisting of all functions $f$ from $\mathcal{O}(U)$  which are stem functions. 

As un example, if $\Delta\subset\C$ is an open disk centered at $0$, each function $F\in\mathcal{O}_s(\Delta,\Kk_n)$ can be represented 
as a convergent series $F(\zeta)=\sum_{k\ge 0}a_k\zeta^k,\,\zeta\in\Delta$, with $a_k\in\Cc_n$ for all $k\ge 0$. 


\begin{Def}\label{vect_fc}\rm Let $U\subset\C$ be a conjugate symmetric open set,  and let $F\in\mathcal{O}(U,\Kk_n)$. For every $\kappa\in U_\sigma$ we set 
\begin{equation}\label{Cauchy_vect}
C[F](\kappa)=\frac{1}{2\pi i}\int_\Gamma F(\zeta)(\zeta-\kappa)^{-1}d\zeta,
\end{equation}
where $\Gamma$ is the boundary of a Cauchy domain in $U$ containing the spectrum $\sigma({\kappa})$. The function $C[F]:U_\sigma\mapsto\Kk_n$
is called the ($\Kk_n$-{\it valued}) {\it Cauchy transform} of the function $F\in\mathcal{O}(U,\Kk_n)$. Clearly, the function $C[F]$ does not depend on the choice of $\Gamma$ because the function $U\setminus\sigma(\kappa)\ni\zeta\mapsto 
F(\zeta)(\zeta-\kappa)^{-1}\in\Kk_n$ is analytic. 

We shall  put
\begin{equation}\label{imagqCt}
\mathcal{R}(U_\sigma,\Kk_n)=\{C[F];F\in\mathcal{O}(U,\Kk_n)\}.
\end{equation}
\end{Def} 


\begin{Pro}\label{right_reg0} Let $U\subset\C$ be open and conjugate symmetric, and let $F\in\mathcal{O}(U,\Kk_n)$. Then  function  $C[F]\in \mathcal{R}(U_\sigma,\Kk_n)$ is slice regular on $U_\sigma$.
\end{Pro}

{\it Proof.}\, Let $F\in\mathcal{O}(U,\Kk_n)$, let
$\kappa\in U_\sigma$ and let  $\Delta\ni\sigma(\kappa)$ be a conjugate symmetric Cauchy domain in $U$, whose boundary is denoted by $\Gamma$. We use the representation of $C[F](\kappa)$
given by (\ref{Cauchy_vect}). Because we have
$$
\bar{\partial}_\mathfrak{s}((\zeta-\kappa)^{-1})=\bar{\partial}_\mathfrak{s}((\zeta-x-y\mathfrak{s})^{-1})=0
$$
for $\kappa=x+y\mathfrak{s}\in\Delta_\sigma\cap(\R+\R\mathfrak{s})$,
via Example \ref{SliceregCk}, we infer that 
$$
\bar{\partial}_\mathfrak{s}(C[F](\kappa))=\frac{1}{2\pi i}\int_\Gamma F(\zeta)\bar{\partial}_\mathfrak{s}((\zeta-\kappa)^{-1})
d\zeta=0,
$$
which implies the assertion.


\begin{Rem}\label{left_mult}\rm
(1)  Because the function $F$ does not necessarily commute with the left multiplication by 
$\mathfrak{s}\in\S_n$, the choice of the right multiplication in the 
slice regularity is necessary to get the stated property of 
$C[F]$. 

(2)  Let $r>0$ and let $U\supset\{\zeta\in\C;\vert\zeta\vert\le r\}$ be a conjugate symmetric open set. Then for every
$F\in\mathcal{O}(U,\Kk_n)$ one has
$$
C[F](\kappa)=\sum_{n\ge0}\frac{F^{(n)}(0)}{n!} \kappa^n,\,\,\vert\kappa \vert<r,
$$
where the series is absolutely convergent. Of course, using the
convergent series $(\zeta-\kappa)^{-1}=
\sum_{n\ge0}\zeta^{-n-1}\kappa^n$ in $\{\zeta;\vert\zeta\vert=r\}$,
the assertion follows easily. Moreover, by Proposition \ref{right_reg0}, the function $C[F]$ is a slice regular $\Kk_n$-valued function in $U_\sigma$. Nevertheless, we are particularly interested in slice regular $\Cc_n$-valued functions.
\end{Rem}

The next result is a version of Theorem 4 from \cite{Vas2}, stated in a quaternionic context. 


\begin{Thm}\label{vect_afc} Let $U\subset\C$ be a conjugate symmetric open set 
and let $F\in\mathcal{O}(U,\Kk_n)$. The Cauchy transform $C[F]$ is 
$\Cc_n$-valued if and only if $F\in\mathcal{O}_s(U,\Kk_n)$.
\end{Thm} 

{\it Proof.\,} We first fix a paravector $\kappa\in U_\sigma\setminus\R$.  If  $\sigma({\kappa})=\{s_+,s_-\}$, the points $s_+,s_-$ are distinct and  not real. We then choose an $r>0$ sufficiently small such that, setting $D_\pm:=\{\zeta\in U;\vert\zeta-s_\pm\vert\le r\}$, we
have 
$D_\pm\subset U$ and $D_+\cap D_-=\emptyset$. Then 
$$
C[F](\kappa)=\frac{1}{2\pi i}\int_{\Gamma_+} F(\zeta)(\zeta-\kappa)^{-1}d\zeta+\frac{1}{2\pi i}\int_{\Gamma_-} F(\zeta)(\zeta-\kappa)^{-1}d\zeta,
$$
where $\Gamma_\pm$ is the boundary of $D_\pm$. We may write 
$F(\zeta)=\sum_{k\ge0}(\zeta-s_+)^ka_k$ with $\zeta\in D_+$,
$a_k\in\Kk_n$ for all $k\ge0$, as a uniformly convergent series.
Similarly, $F(\zeta)=\sum_{k\ge0}(\zeta-s_-)^kb_k$ with $\zeta\in D_-$, $b_k\in\Kk_n$ for all $k\ge0$, as a uniformly convergent series. 
 
Note that
$$
\frac{1}{2\pi i}\int_{\Gamma_+} F(\zeta)(\zeta-\kappa)^{-1}d\zeta=\sum_{k\ge0}\left(a_k\frac{1}{2\pi i}\int_{\Gamma_+} (\zeta-s_+)^k(\zeta-\kappa)^{-1}d\zeta\right)=
a_0\iota_+(\mathfrak{s}_{\tilde{\kappa}})
$$
because  we have
$$
\frac{1}{2\pi i}\int_{\Gamma_+} (\zeta-s_+)^k(\zeta-\kappa)^{-1}d\zeta=(\kappa-s_+)^k\iota_+(\mathfrak{s}_{\tilde{\kappa}})
$$
by the analytic functional calculus of $\kappa$ (see also
Lemma \ref{spproj}), which is equal
to $\iota_+(\mathfrak{s}_{\tilde{\kappa}})$ when $k=0$,  and it is equal to $0$ when $k\ge 1$, via the equality ${\kappa}\iota_+(\mathfrak{s}_{\tilde{\kappa}})=s_+\iota_+(\mathfrak{s}_{\tilde{\kappa}})$ 

\medskip

Similarly
$$
\frac{1}{2\pi i}\int_{\Gamma_-} F(\zeta)(\zeta-\kappa)^{-1}d\zeta=\sum_{k\ge0}\left(b_k\frac{1}{2\pi i}\int_{\Gamma_-} (\zeta-s_-)^k(\zeta-\kappa)^{-1}d\zeta\right)=b_0\iota_-(\mathfrak{s}_{\tilde{\kappa}})
$$
because, as above, we have
$$
\frac{1}{2\pi i}\int_{\Gamma_-} (\zeta-s_-)^k
(\zeta-\kappa)^{-1})d\zeta=(\kappa-s_-)^k\iota_-(\mathfrak{s}_{\tilde{\kappa}}),
$$
which is equal $\iota_-(\mathfrak{s}_{\tilde{\kappa}})$ when $k=0$, and it is equal to $0$ when $k\ge 1$. Consequently,
$$
C[F](\kappa)=F(s_+)\iota_+(\mathfrak{s}_{\tilde{\kappa}})+F(s_-)\iota_-(\mathfrak{s}_{\tilde{\kappa}}),
$$
and the right hand side of this equality coincides with the expression from formula (\ref{GFC2}). 
 
Assume now that
$\sigma(\kappa)=\{s\}$,  where  $s:=s_+=s_-\in\R$. Then necessarily 
$\kappa=s$, so fixing  an $r>0$ such that the set 
$D:=\{\zeta\in U;\vert\zeta-s\vert\le r\}\subset U$, 
whose boundary is denoted by $\Gamma$, we have
$$
C[F](\kappa)=\frac{1}{2\pi i}\int_{\Gamma} F(\zeta)(\zeta-\kappa)^{-1}d\zeta=F(s),
$$
via the usual analytic functional calculus.

 In all of these situations, the element $C[F](\kappa)$ is equal to  the right hand side of formula (\ref{GFC2}). Therefore, we must have 
$C[F](\kappa)\in\Cc_n$ if and only if $F(s_+)=\overline{F(s_-)}$, via Theorem \ref{sym_spec}. As every point $\lambda\in U$ is the eigenvalue of a certain paravector in $U_\sigma$, we deduce that 
$C[F](\kappa)\in\Cc_n$ for all $\kappa\in U_\sigma$ if and only if $F:U\mapsto \Kk_n$ is a stem function.  
 

\begin{Rem}\label{bi_not}\rm (1) It follows from the proof of the previous theorem 
that the element $C[F](\kappa))$, given by formula (\ref{Cauchy_vect}), coincides with the element $F_\sigma(\kappa)$ given by (\ref{GFC2}). To unify the notation, from now on this
element will be denoted by $F_\sigma(\kappa)$, whenever $F$ is a stem function analytic or not.

(2) An important particular case is when let $f:U\mapsto\C$ is an analytic function, where  $U\subset\C$ is a conjugate symmetric open set. In this case we may also consider the  Cauchy transform of $f$ given by
 
\begin{equation}\label{Cauchy_vect1}
C[f](\kappa)=\frac{1}{2\pi i}\int_\Gamma f(\zeta)(\zeta-\kappa)^{-1}d\zeta,
\end{equation}
where  $\Gamma$ is the boundary of a Cauchy domain in $U$ containing the spectrum $\sigma(\kappa)$. According to Theorem 
\ref{vect_fc},  we have  $C[f](\kappa)\in\Cc_n$ if and only if $f$ is a stem function, that is 
$f\in\mathcal{O}_s(U)$. Of course, in this case 
 we may (and shall) also use the notation $C[f]=f_\sigma$, and we have, in fact, $f_\sigma\in\mathcal{IF}(U_\sigma,\Cc_n)$, by
Corollary \ref{intrinsic0}. 
\end{Rem}


\section{Analytic Functional Calculus for Stem\\  Functions}
\label{AFCQ}

Let $\Omega\subset\Pp_n$ be a spectrally saturated open set, and 
let $U=\mathfrak{S}(\Omega)\subset\C$ (which is conjugate symmetric, and also open, by
Remark \ref{consymr}(2)). We introduce the notation
$$
\mathcal{R}_{s,n}(\Omega)=\{f_\sigma; f\in\mathcal{O}_s(U)\},
$$

$$
\mathcal{R}_s(\Omega,\Cc_n)=\{F_\sigma; F\in\mathcal{O}_s(U,\Kk_n)\},
$$
which are $\R$-vector spaces.

In fact, these $\R$-vector spaces have some important 
properties, as already noticed in a quaternionic version of the next theorem (see Theorem 5 in \cite{Vas2}).


\begin{Thm}\label{C_afc} Let $\Omega\subset\Pp_n$ be a spectrally
saturated open set, and let $U=\mathfrak{S}(\Omega)$.
The space $\mathcal{R}_{s,n}(\Omega)$ is a unital commutative $\R$-algebra,  the space $\mathcal{R}_s(\Omega,\Cc_n)$ is a right
$\mathcal{R}_{s,n}(\Omega)$-module, the linear map
$$
{\mathcal O}_s(U,\Kk_n)\ni F\mapsto F_\sigma\in\mathcal{R}_s(\Omega,\Cc_n)
$$ 
is a right module isomorphism, and  its restriction
$$
{\mathcal O}_s(U)\ni f\mapsto f_\sigma\in\mathcal{R}_{s,n}(\Omega)
$$ 
is an $\R$-algebra isomorphism.

Moreover, for every polynomial\,\,$P(\zeta)=\sum_{n=0}^m a_n\zeta^n,\,\zeta\in\C$, with $a_n\in\Cc_n$  for all $n=0,1,\ldots,m$, we have  $P_\sigma(\kappa)=\sum_{n=0}^m a_n \kappa^n\in\Cc_n$ for all $\kappa\in\Pp_n$.   
\end{Thm}

{\it Proof.\,} Thanks to Theorem \ref{vect_afc}, this statement 
is a particular case of Theorem \ref{gfc}. Indeed, the 
 $\R$-linear maps
$$
{\mathcal O}_s(U,\Kk_n)\ni F\mapsto F_\sigma\in\mathcal{R}_s(\Omega,\Cc_n),\,
{\mathcal O}_s(U)\ni f\mapsto f_\sigma\in\mathcal{R}_{s,n}(\Omega),\,\,
$$
are restrictions of the maps
$$
{\mathcal S}(U,\Kk_n)\ni F\mapsto F_\sigma\in\mathcal{F}(\Omega,\Cc_n),\,\,
{\mathcal S}(U)\ni f\mapsto f_\sigma\in\mathcal{IF}(\Omega,\Cc_n),
$$ 
respectively. Moreover, they are $\R$-isomorphisms, the latter being actually unital and  multiplicative.   
Note that, in particular,  for every polynomial $P(\zeta)=\sum_{n=0}^m a_n\zeta^n$ with $a_n\in\Cc_n$  for all $n=0,1,\ldots,m$, we have  $P_\sigma(\kappa)=\sum_{n=0}^m a_n \kappa^n\in\Cc_n$ for all 
$\kappa\in\Pp_n$. 


\begin{Rem}\label{n-deriv}\rm For every function $F\in\mathcal{O}_s(U,\Kk_n)$, the derivatives $F^{(n)}$ also belong to 
$\mathcal{O}_s(U,\Kk_n)$, where $U\subset\C$ is a conjugate symmetric open set.  

Now fixing $F\in\mathcal{O}_s(U,\Kk_n)$, we may define its 
{\it extended derivatives} with respect to the paravector variable via the formula

\begin{equation}\label{Cauchy_vect_deriv}
F^{(n)}_\sigma(\kappa)=\frac{1}{2\pi i}\int_\Gamma F^{(n)}(\zeta)(\zeta-\kappa)^{-1}d\zeta,
\end{equation}
for the boundary $\Gamma$ of a Cauchy domain $\Delta\subset U$, $n\ge0$ an arbitrary integer, and $\sigma(\kappa)\subset\Delta$. 

In particular, if $\Delta$ is a disk centered at zero and 
$F\in\mathcal{O}_s(\Delta,\Kk_n)$, so we have a representation of
$F$   as a convergent series $\sum_{m\ge0}a_k\zeta^k$ with coefficients in $\Cc_n$, then 
(\ref{Cauchy_vect_deriv}) gives the equality
$F'_\sigma(\kappa)=\sum_{m\ge1}ma_m\kappa^{m-1}$, which looks like a (formal) derivative of the function $F_\sigma(\kappa)=\sum_{m\ge0}a_k\kappa^{m}$. 
\end{Rem}


\begin{Rem}\label{deriv}\rm As already noticed in the framework of 
\cite{Vas2}, Theorem \ref{C_afc} suggests 
a definition for $\Cc_n$-valued "analytic functions`` as elements of the set $\mathcal{R}_s(\Omega,\Cc_n)$, where $\Omega$ is a 
spectrally saturated open subset of $\Pp_n$. 
Because the expression "analytic function`` is quite improper 
in this context, the elements of $\mathcal{R}_s(\Omega,\Cc_n)$ will be called {\it $\Cc_n$-regular functions} on $\Omega$. As
shown by Theorem \ref{C_afc}, 
the functions from $\mathcal{R}_s(\Omega,\Cc_n)$ are the  Cauchy transforms of the stem functions from $\mathcal{O}_s(U,\Kk_n)$, with $U=\mathfrak{S}(\Omega)$. 

Except for Theorem \ref{C_afc}, many properties of $\Cc_n$-regular 
functions can be obtained directly from the definition, by 
recapturing the corresponding results from \cite{Vas2}. We omit the details.
\end{Rem}


\begin{Rem}\label{repres_form}\rm Let $U\subset\C$ be a conjugate symmetric open set, let $x,y\in\R$ with $y\neq0$ and $z_\pm=x\pm iy\in U$, let  $F\in\mathcal{O}_s(U,\Kk_n)$, and let $\mathfrak{s}\in\S_n$. We shall apply some arguments similar to
those from Corollary \ref{intrinsic0}, in a non commuting context. 

Assuming $y>0$, we consider the paravectors 
$\kappa_\pm=x\pm y \mathfrak{s}$ for which
$s_+(\kappa_\pm)=x+iy,\,s_-(\kappa_\pm)=x-iy$. As we have
$\tilde{\kappa}_\pm=\pm y\mathfrak{s}$, then  
$s_{\tilde{\kappa}_\pm}=\pm\mathfrak{s}$, and $\iota_\pm(s_{\tilde{\kappa}_+})=(1\mp i\mathfrak{s})/2, 
\iota_\pm(s_{\tilde{\kappa}_-})=(1\pm i\mathfrak{s})/2.$ Therefore,
$$
2F_\sigma({\kappa}_+)=F(z_+)(1-i\mathfrak{s})+F(z_-)(1+i\mathfrak{s}),
$$
$$
2F_\sigma({\kappa}_-)=F(z_+)((1+i\mathfrak{s})+F(z_-)(1-i\mathfrak{s}).
$$
From these equations we deduce that
\begin{equation}\label{reprez1}
2F(z_+)=F_\sigma({\kappa}_+)(1-i\mathfrak{s})+F_\sigma({\kappa}_-)(1+i\mathfrak{s}),
\end{equation}
\begin{equation}\label{reprez2}
2F(z_-)=F_\sigma({\kappa}_+)(1+i\mathfrak{s})+F_\sigma({\kappa}_-)(1-i\mathfrak{s}).
\end{equation}

If $y<0$, for the paravectors ${\kappa}_\pm=x\pm y \mathfrak{s}$
we have $s_+({\kappa}_\pm)=x-iy,\,s_-({\kappa}_\pm)=x+iy$.
Moreover, as $\tilde{\kappa}_\pm=\pm y\mathfrak{s}$, then  
$s_{\tilde{\kappa}_\pm}=\mp\mathfrak{s}$, and $\iota_\pm(s_{\tilde{\kappa}_+})=(1\pm i\mathfrak{s})/2,\break 
\iota_\pm(s_{\tilde{\kappa}_-})=(1\mp i\mathfrak{s})/2.$
Therefore
$$
2F_\sigma({\kappa}_+)=F(z_-)(1+i\mathfrak{s})+F(z_+)(1-i\mathfrak{s}),
$$
$$
2F_\sigma({\kappa}_-)=F(z_-)((1-i\mathfrak{s})+F(z_+)(1+i\mathfrak{s}).
$$
These formulas lead again to equations (\ref{reprez1}) and 
(\ref{reprez2}). Consequently, we have the following (see also \cite{CoSaSt}, Theorem 2.2.18 for a similar result, and 
\cite{Vas2} for the corresponding result in the quaternionic context).

The next proposition, and  Remark \ref{repres_form} as well, have their counterparts in \cite{Vas2}, stated as  Proposition 3 and
Remark 10, respectively. For a similar result see also Theorem 2.2.18 from \cite{CoSaSt}.


\begin{Pro} Let $U\subset\C$ be a conjugate symmetric open set, let $x,y\in\R$ with $x\pm iy\in U$, let $\mathfrak{s}\in\S_n$, and let $F\in\mathcal{O}_s(U,\Kk_n)$. Then we have the formulas
\begin{equation}\label{reprez}
F(x\pm iy)=F_\sigma(x\pm y\mathfrak{s})\left(\frac{1\mp i\mathfrak{s}}{2}\right)+
F_\sigma(x\mp y\mathfrak{s})\left(\frac{1\pm i\mathfrak{s}}{2}\right).
\end{equation}
\end{Pro} 

As the proof has been previously done, we only note that 
equality (\ref{reprez}) also holds  for $y=0$.
\end{Rem}


\begin{Lem}\label{equiv-ons-dom} Let $U\subset\C$ be a conjugate symmetric open set, let $\mathfrak{s}\in\S_n$ be fixed,
and let $\Psi:U_\mathfrak{s}\mapsto\Cc_n$ be  such that $\bar{\partial}_{\pm\mathfrak{s}}\Psi=0$. Then there exists  a function $\Phi\in\mathcal{R}_s(U_\sigma,\Cc_n)$ with $\Psi=\Phi\vert U_\mathfrak{s}$,
where  $ U_\mathfrak{s}=\{x+y\mathfrak{s}; x+iy\in U\}$. 
\end{Lem}

{\it Proof.}\, For arbitrary points $z_\pm=x\pm iy\in U$  with $x,y(\neq0)\in 
\R$, as in Remark \ref{repres_form}, we consider the paravectors ${\kappa}_\pm=x\pm y \mathfrak{s}$, so $s_+( {\kappa}_\pm)=x+i\vert y\vert$, and $s_-( {\kappa}_\pm)=x-i\vert y\vert$. Inspired by formula (\ref{reprez}),  we set
$$
2F(z_+)=\Psi({\kappa}_+)(1-i\mathfrak{s})+\Psi({\kappa}_-)(1+i\mathfrak{s}),
$$
$$
2F(z_-)=\Psi({\kappa}_+)(1+i\mathfrak{s})+\Psi({\kappa}_-)(1-i\mathfrak{s}).
$$
 
Then we have
$$
2\frac{\partial F(z_+)}{\partial x}=\frac{\partial\Psi({\kappa}_+) }{\partial x}(1-i\mathfrak{s})+\frac{\partial\Psi({\kappa}_-) }
{\partial x}(1+i\mathfrak{s}),
$$
and
$$
2i\frac{\partial F(z_+)}{\partial y}=\frac{\partial\Psi({\kappa}_+) }{\partial y}\mathfrak{s}(1-i\mathfrak{s})+\frac{\partial\Psi({\kappa}_-)}{\partial y}(-\mathfrak{s})(1+i\mathfrak{s}),
$$
because $i(1-i\mathfrak{s})=\mathfrak{s}(1-i\mathfrak{s})$
and $i(1+i\mathfrak{s})=-\mathfrak{s}(1+i\mathfrak{s})$.

Therefore,
$$
\frac{\partial F(z_+)}{\partial x}+i\frac{\partial F(z_+)}{\partial y}=(\bar{\partial}_\mathfrak{s}\Psi({\kappa}_+))(1-i\mathfrak{s})+(\bar{\partial}_\mathfrak{-s}\Psi({\kappa}_-))(1+i\mathfrak{s})=0,
$$
showing that the function $z_+\mapsto F(z_+)$ is analytic in $U$. 

Because $\overline{F(z_-)}=\overline{F(\overline{z_+})}=F(z_+)$, and when $y=0$ we have
$\overline{F(z_-)}=F(z_+)=F(x)$, we have constructed a function
$F\in\mathcal{O}_s(U,\Kk_n)$. Hence, taking $\Phi=F_\sigma$, we have
$\Phi\in\mathcal{R}_s(U_\sigma,\Cc_n)$ with $\Psi=\Phi\vert U_\mathfrak{s}$, via Remark \ref{repres_form}.
\medskip


The next theorem is a version of  Theorem 6 from \cite{Vas2}.
 
\begin{Thm}\label{equiv-ons-dom1} Let $\Omega\subset\Pp_n$ be a spectrally saturated open set, and let $\Phi:\Omega\mapsto\Cc_n$.
 The following conditions are equivalent:

$(i)$ $\Phi$ is  a slice regular function;  

$(ii)$ $\Phi\in\mathcal{R}_s(\Omega,\Cc_n)$, that is, $\Phi$ is 
$\Cc_n$-regular. 
\end{Thm}

{\it Proof.}\, If $\Phi\in\mathcal{R}_s(\Omega,\Cc_n)$, then $\Phi$
is slice regular, by Proposition \ref{right_reg0}, so $(ii)\Rightarrow(i)$.
 
Conversely, let $\Phi$ be slice regular in $\Omega$. Fixing an 
$\mathfrak{s}\in\S_n$, we have $\bar{\partial}_{\pm\mathfrak{s}}\Phi_{\mathfrak{s}}=0$, where $\Phi_{\mathfrak{s}}=\Phi\vert U_{\mathfrak{s}}$. It follows from Lemma \ref{equiv-ons-dom}
that there exists $\Psi\in\mathcal{R}(U_\sigma,\Cc_n)$ with
$\Psi_{\mathfrak{s}}=\Phi_{\mathfrak{s}}$. This implies that 
$\Phi=\Psi$, because both $\Phi,\Psi$ are uniquely determined by 
$\Phi_{\mathfrak{s}}, \Psi_{\mathfrak{s}}$, respectively, the former by (the right hand version of) Lemma 2.2.24  in \cite{CoSaSt}, and the latter by  Remark \ref{zeros}. Consequently, we also have $(i)\Rightarrow(ii)$.
\medskip

\noindent{\bf Remark}  A concept of ''Cliffordian holomorphic function`` also appears in  \cite{LaRa}, in a different context. 


\section{Spectrum of Clifford  Operators}

Let $\V$ be a Banach $Cl$-space, and let $\V_\C=\V+i\V$ its complexification, endoved with the norm $\Vert x+iy\Vert=\Vert x\Vert+\Vert y\Vert$, for all $x,y\in\V$, where  $\Vert*\Vert$ is the norm of $\V$.
We denote by $C$ the {\it conjugation} on $\V_\C$, that is, the map $C(x+iy)=x-iy$ for all $x,y\in\V$, which is an $\R$-linear map whose square is the identity. 

 As in Subsection 2.3,
for every $T\in \B(\V)$, we consider its natural ''complex  extension`` to $\V_\C$ given by  $T_\C(x+iy)=Tx+iTy$, for all $x,y\in\V$,  which is at least $\C$-linear (if $T\in\B^{\rm r}(\V)$,  $T_\C$  is  right $\Kk_n$-linear),  so $T_\C\in\B(\V_\C)$.  As already noticed, the map $\B(\V)\ni T\mapsto T_\C\in\B(\V_\C)$ is a unital injective morphism of real algebras.

Assuming that $\V$ is a  Banach $Cl$-space implies that
$\mathcal{B}^{\rm r}(\V)$ is a unital real Banach $Cl$-algebra 
(that is, a Banach algebra which also a Banach $Cl$-space), via the 
algebraic operations  $(\a T)(x)= \a T(x)$, and $(T\a)(x)=T(\a x)$ for all
$\a\in\Cc_n$ and $x\in\V$. The complexification $\mathcal{B}^{\rm r}(\V)_\C$ of $\mathcal{B}^{\rm r}(\V)$
 is, in particular, a unital complex Banach algebra, with the product
\smallskip

$(T_1+iT_2)(S_1+iS_2)=T_1S_1-T_2S_2+i(T_1S_2+T_2S_1),\,\,T_1,T_2,S_1,S_2\in\mathcal{B^{\rm r}(V)},$ 

\noindent and  a fixed norm, say   $\Vert (T_1+iT_2)\Vert=\Vert T_1\Vert+\Vert T_2\Vert,
T_1,T_2\in\B^{\rm r}(\V)$. 

 Also note that the  complex numbers, regarded as elements of $\mathcal{B}^{\rm r}(\V)_\C$, commute with
the elements of $\mathcal{B}^{\rm r}(\V)$.


\begin{Rem}\label{iso0}\rm  For every  $S\in\B(\V_\C)$ we put $S^\flat=CSC\in\B(\V_\C)$. The assignment $S\mapsto S^\flat$ is a conjugate linear automorphism of the algebra  $\B(\V_\C)$, whose square is the identity operator. In fact, the map $S\mapsto S^\flat$ is a conjugation of $\B(\V)$, induced by $C$. Moreover, $S^\flat=S$ if and only if $S(\V)\subset\V$.  In particular, we have $S=S_1+iS_2$ with $ S_j(\V)\subset\V,\,j=1,2,$ uniquely determined. Its action on the space $\V_\C$ is given by
$S(x+iy)=S_1x-S_2y+i(S_1y+S_2x)$ for all $x,y\in\V$.

Because  $CR_\a=R_\a C$ for all $\a\in\Pp_n$, it follows that if  $S\in\mathcal{B}^{\rm r}(\V_\C)$, then $S^\flat\in\B^{\rm r}(\V_\C)$. In particular, we have $(S+S^\flat)(\V)\subset\V$,  $i(S-S^\flat)(\V)\subset\V$, and
$(T_\C)^\flat=T_\C$ for all $T\in\B^{\rm r}(\V)$.  Note also that the map
$$
\B^{\rm r}(\V)\ni T\mapsto T_\C\in\{S=S^\flat;S\in \B^{\rm r}(\V_\C)\}
$$
is actually a real unital  algebra  isomorphism, since its surjectivity follows from the equality $(S\vert  \V)_\C=S$ whenever $S=S^\flat\in \B^{\rm r}(\V_\C)$. This implies that  the algebras $\B^{\rm r}(\V_\C)$ and $\B^{\rm r}(\V)_\C$ are isomorphic. 
 This isomorphism is given by the assignment
$$
\B^{\rm r}(\V)_\C\ni T_1+iT_2\mapsto T_{1\C}+iT_{2\C}\in\B^{\rm r}(\V_\C)
$$
which is  is actually an algebra isomorphism, via a direct calculation.

The continuity of this assignment  is also clear, and therefore it is a Banach algebra isomorphism.
For this reason, may  identify the algebras  $\B^{\rm r}(\V_\C)$ and $\B^{\rm r}(\V)_\C$.
As already noticed above, the real algebras  $\B^{\rm r}(\V)$ and $\{S\in \B^{\rm r}(\V_\C);S=S^\flat\}$ may and will  be also
identified.
\end{Rem}

The operators from the algebra $\B^{\rm r}(\V)$ will be sometimes called {\it Clifford operators}, or simply $Cl$-{\it operators}.

Looking at  Definition 3.1.4  from \cite{CoSaSt}, we can give  the folowing. 


\begin{Def}\label{C-spectrum}\rm For a given operator $T\in B^{\rm r}(\V)$, the set 
$$
\sigma_{Cl}(T):=\{\kappa\in\Pp_n; T^2-2\Re(\kappa)T+\vert\kappa\vert^2)\,\,
{\rm not}\,\,{\rm invertible}\}
$$
is called the {\it Clifford} (or $Cl$-){\it spectrum} of $T$.
 
The complement $\rho_C(T)=\Pp_n\setminus\sigma_{Cl}(T)$ is called 
the   {\it Clifford} (or $Cl$-){\it  resolvent of $T$}.
\end{Def}

 Note that, if $\a\in\sigma_{Cl}(T)$, then 
 $\{\b\in\Pp_n;\sigma(\b)=\sigma(\a)\}\subset\sigma_{Cl}(T)$.
  
Since every operator $T\in\mathcal{B}^{\rm r}(\V)$ is, in particular,  $\R$-linear, we also have a {\it complex resolvent},
defined by
$$
\rho_\C(T)=\{\lambda\in\C;(T^2-2\Re(\lambda)T+\vert\lambda\vert^2)^{-1}\in\mathcal{B}^{\rm r}(\V)= 
$$
$$
\{\lambda\in\C;(\lambda-T_\C)^{-1}\in \mathcal{B}^{\rm r}(\V_\C)\}=\rho(T_\C),
$$  
and the associated {\it complex spectrum} $\sigma_\C(T)=\sigma(T_\C)$ as well.

Note that both sets $\sigma_\C(T)$ and $\rho_\C(T)$ are conjugate symmetric.
 
There exists a strong connexion between $\sigma_{Cl}(T)$ and 
$\sigma_\C(T)$. In fact, the set $\sigma_\C(T)$ looks like a ''complex border`` of the set $\sigma_{Cl}(T)$. 
Specifically, we can prove the following.

\begin{Lem}\label{spec_eg0} For every $T\in\mathcal{B}^{\rm r}(\V)$ we have the
equalities 
\begin{equation}\label{spec_eg} 
\sigma_{Cl}(T)=\{\a\in\Pp_n;\sigma_\C(T)\cap\sigma(\a)\neq\emptyset\}.
\end{equation}
and 
\begin{equation}\label{spec_eg1}
\sigma_\C(T)=\{\lambda\in\sigma(\a);\a\in\sigma_{Cl}(T)\}.
\end{equation}
 \end{Lem}

{\it Proof.} Let us prove (\ref{spec_eg}). If $\kappa\in \sigma_{Cl}(T)$, and so the $T^2-2\Re(\kappa)T+\vert\kappa\vert^2$ is not invertible, choosing  $\lambda\in\{\Re(\kappa)\pm i\vert\Im(\kappa)\vert\}=\sigma(\kappa)$, we clearly have $T^2-2\Re(\lambda)T+\vert\lambda\vert^2$ not invertible, implying $\lambda\in\sigma_\C(T)\cap\sigma(\kappa)\neq\emptyset$. 

Conversely, if for some $\kappa\in\Pp_n$ there exists $\lambda\in\sigma_\C(T)\cap\sigma(\kappa)$, and so $T^2-2\Re(\lambda)T+\vert\lambda\vert^2=T^2-2\Re(\kappa)T+\vert\kappa\vert^2$ is not invertible, we must have $\kappa\in\sigma_{Cl}(T)$. 

 We now prove (\ref{spec_eg1}). Let $\lambda\in\sigma_\C(T)$, so the operator $T^2-2(\Re\lambda)T+\vert\lambda\vert^2$ is not invertible. Setting $\kappa=\Re(\lambda)+\vert\Im\lambda\vert\mathfrak{s}$, with $\mathfrak{s}\in\mathbb{S}_n$, we have 
$\lambda\in\sigma(\kappa)$. Moreover, $T^2+2\Re(\kappa)T+\vert\kappa\vert^2$ is not 
invertible, and so $\kappa\in\sigma_ {Cl}(T)$. 

Conversely, if $\lambda\in\sigma(\kappa)$ for some $\kappa\in\sigma_{Cl}(T)$, then
$\lambda\in\{\Re(\kappa)\pm i\vert\Im(\kappa)\vert\}$, showing that
$T^2-2\Re(\lambda)T+\vert\lambda\vert^2=T^2+2\Re(\kappa)T+\vert\kappa\vert^2$ is
not invertible. 
\medskip


\begin{Rem}\rm As expected, the set $\sigma_{Cl}(T)$ is nonempty and bounded, which follows from Lemma \ref{spec_eg0}.  In fact, we have the equality
$$
\sigma_{Cl}(T)=\{\Re(\lambda)+\vert\Im(\lambda)\vert\mathfrak{s};\lambda\in\sigma_\C(T),\mathfrak{s}\in\mathbb{S}_n\}.
$$
It is also closed, as a consequence of Definition \ref{C-spectrum}, because the set of invertible elements in $\mathcal{B}^{\rm r}(\V)$ is open.
\end{Rem}

Note that the subset $\sigma_{Cl}(T)$ spectrally saturated (see Definition \ref{consym}(2)).


\section{Analytic Functional Calculus for  Clifford\\ Operators}

Having a concept of spectrum for Clifford operators, an important step for further development is the construction of an analytic functional calculus. We follow the main ideas from \cite{Vas3}.

If $\V$ is a  Banach $Cl$-space, and so each operator $T\in\mathcal{B}^{\rm r}(\V)$ has a complex spectrum $\sigma_\C(T)$,  one can  use the classical Riesz-Dunford functional calculus, in a slightly generalized form (that is, replacing the scalar-valued analytic functions by operator-valued analytic ones, which is a well known idea).


\begin{Rem}\label{afcro}\rm If $\mathcal{V}$ is  Banach $Cl$-space, and 
$T\in \B(\V)$, we have the usual analytic functional calculus for the operator $T_\C\in\mathcal{B}(\mathcal{V}_\C)$
(see \cite{DuSc}). That is, in a slightly generalized form, and for later use, if $U\supset\sigma(T_\C)$ is an open set in $\C$ and 
$F:U\mapsto B(\V_\C)$ is analytic, the (left)  Riesz-Dunford analytic functional calculus is given by the formula
$$
F(T_\C)=\frac{1}{2\pi i}\int_\Gamma F(\zeta)(\zeta-T_\C)^{-1}
d\zeta,
$$  
where $\Gamma$ is the boundary of a Cauchy domain $\Delta$ containing  
$\sigma(T_\C)$ in $U$. In fact, since $\sigma(T_\C)$ is conjugate symmetric, we may and shall assume 
that both $U$ and $\Gamma$ are conjugate symmetric. Because the
function $\zeta\mapsto F(\zeta)(\zeta-T_\C)^{-1}$ is analytic
in $U\setminus\sigma(T_\C)$, the integral does not depend on the
particular choice of the Cauchy domain $\Delta$ containing 
$\sigma(T_\C)$.

A natural question is to find an appropriate condition to  have  $F(T_\C)^\flat=F(T_\C)$, which would
imply the invariance of $\mathcal{V}$ under $F(T_\C)$. 
\end{Rem}


\begin{Rem}\label{stem_anal}\rm If $\A$ is a unital real Banach algebra, 
$\A_\C$ its complexification, and $U\subset\C$ is open, we denote by 
$\mathcal{O}(U,\A_\C)$ the algebra of all analytic $\A_\C$-valued functions. If $U$ is conjugate symmetric, and 
$\A_\C\ni a\mapsto \bar{a}\in
\A_\C$ is its natural conjugation, we denote by $\mathcal{O}_s(U,\A_\C)$ the real subalgebra of $\mathcal{O}(U,\A_\C)$ consisting of those functions $F$ with the property $F(\bar{\zeta})=\overline{F(\zeta)}$
for all $\zeta\in U$.   As in Definition \ref{stem},  such 
functions will   be called ($\A_\C$-{\it valued $)$ stem functions}.

When $\A=\R$, so $\A_\C=\C$, the space $\mathcal{O}_s(U,\C)$ will be denoted by 
$\mathcal{O}_s(U)$, which is a real algebra. Note that 
$\mathcal{O}_s(U,\A_\C)$ is also a two-sided $\mathcal{O}_s(U)$-module. 
\end{Rem}
 
With the notation of Remark \ref{afcro}, we state and prove the following adapted version of Theorem 1 from \cite{Vas3}.


\begin{Thm}\label{afcro1} Let $U\subset\C$ be open and conjugate symmetric. If $F\in\mathcal{O}_s(U,\B(\V_\C))$,
we have $F(T_\C)^\flat=F(T_\C)$ for all $T\in\B(\V)$ with $\sigma_\C(T)\subset U$.

 Moreover, if $F\in\mathcal{O}_s(U,\B^{\rm r}(\V_\C))$, and $T\in\B^{\rm r}(\V)$, then  $F(T_\C)\in \B^{\rm r}(\V_\C)$.
\end{Thm}

{\it Proof.}\, We use the  notation from  Remark \ref{afcro}, ssuming $\Gamma$  conjugate symmetric.  We put 
$\Gamma_\pm:=\Gamma\cap\C_\pm$, where $\C_+$ (resp. $\C_-$) equals to $\{\lambda\in\C;\Im\lambda\ge0\}$ (resp. $\{\lambda\in\C;\Im\lambda\le0\}$). We write $\Gamma_+=\cup_{j=1}^m\Gamma_{j+}$, where $\Gamma_{j+}$ are the connected components
of $\Gamma_+$. Similarly, we write  $\Gamma_-=\cup_{j=1}^m\Gamma_{j-}$, where $\Gamma_{j-}$ are the connected components
of $\Gamma_-$, and $\Gamma_{j-}$ is the reflexion of 
$\Gamma_{j+}$ with respect of the real axis. 

As $\Gamma$ is a finite union of Jordan piecewise smooth closed curves, 
for each index
$j$ we have a parametrization $\phi_j:[0,1]\mapsto\C$, positively oriented, such that
$\phi_j([0,1])=\Gamma_{j+}$. Taking into account that the function $t\mapsto\overline{\phi_j(t)}$ is a 
parametrization of $\Gamma_{j-}$ negatively oriented, and setting 
$\Gamma_j=\Gamma_{j+}\cup\Gamma_{j-}$, we can write
$$
F_j(T_\C):=\frac{1}{2\pi i}\int_{\Gamma_j} F(\zeta)(\zeta-T_\C)^{-1}
d\zeta=
$$
$$
\frac{1}{2\pi i}\int_0^1 F(\phi_j(t))(\phi_j(t)-T_\C)^{-1}
\phi_j'(t)dt
$$
$$
-\frac{1}{2\pi i}\int_0^1 F(\overline{\phi_j(t)})(\overline{\phi_j(t)}-T_\C)^{-1}\overline{\phi_j'(t)}dt.
$$
Therefore,
$$
F_j(T_\C)^\flat=
-\frac{1}{2\pi i}\int_0^1 F(\phi_j(t))^\flat(\overline{\phi_j(t)}-T_\C)^{-1}\overline{\phi_j'(t)}dt
$$
$$
+\frac{1}{2\pi i}\int_0^1 F(\overline{\phi_j(t)})^\flat(\phi_j(t)-T_\C)^{-1}\phi_j'(t)dt.
$$
According to our assumption on the function $F$, we obtain
$F_j(T_\C)=F_j(T_\C)^\flat$ for all $j$, and therefore 
$$
F(T_\C)^\flat=\sum_{j=1}^mF_j(T_\C)^\flat=\sum_{j=1}^mF_j(T_\C)=F(T_\C). $$

Moreover, if  $T\in\B^{\rm r}(\V)$,  because $R_\a  F(\zeta)(\zeta-T_\C)^{-1}= F(\zeta)(\zeta-T_\C)^{-1}R_\a$ for all $\a\in\Pp_n$, via the definition of $\B^{\rm r}(\V_\C)$, it follows that $F(T_\C)\in \B^{\rm r}(\V_\C)$.
\medskip

In the next result, we identify the algebra  $\B^{\bf r}(\V)$ with a real subalgebra of 
$\B^{\bf r}(\V)_\C$, in turn identified with $\B^{\bf r}(\V_\C)$  (see Remark \ref{iso0}). In this case, when $F\in \mathcal{O}_s(U,\B^{\bf r}(\V)_\C)$, we shall write
$$
F(T)=\frac{1}{2\pi i}\int_\Gamma F(\zeta)(\zeta-T)^{-1}
d\zeta,
$$  
noting that the right hand side of this formula belongs to 
$\B^{\bf r}(\V)$, as a consequence of Theorem \ref{afcro1}.

The following result expresses the {\it (left)  analytic functional calculus } of a given operator from $ \B^{\rm r}(\V)$  with $ \B^{\rm r}(\V)_\C$-valued stem functions. It is a version of Theorem 4 from \cite{Vas3}, proved  in a quaternionic context.


\begin{Thm}\label{R_afc} Let $\V$ be a Banach $Cl$-space, let $U\subset\C$ be a conjugate symmetric open set, and let 
$T\in \B^{\rm r}(\V)$, with $\sigma_\C(T)\subset U$.  Then the assignment 
$$
{\mathcal O}_s(U,   \B^{\rm r}(\V)_\C)\ni F\mapsto F(T)\in   \B^{\rm r}(\V)
$$ 
is an $\R$-linear map, and the map
$$
{\mathcal O}_s(U)\ni f\mapsto f(T)\in   \B^{\rm r}(\V)
$$
is a unital real algebra morphism.

Moreover, the following properties hold true:

(1) for all   $F\in\mathcal{O}_s(U,  \B^{\bf r}(\mathcal{V})_\C),\, f\in{\mathcal O}_s(U)$, we have $(Ff)(T)=F(T)f(T)$.

(2) for every polynomial $P(\zeta)=\sum_{n=0}^m A_n\zeta^n,\,\zeta\in\C$, with $A_n\in \B^{\bf r}(\mathcal{V})$  for all $n=0,1,\ldots,m$, we have  $P(T)=\sum_{n=0}^m A_n T^n\in  \B^{\bf r}(\mathcal{V})$.  
\end{Thm}

{\it Proof.\,} The arguments are more or less standard (see
\cite{DuSc}). The $\R$-linearity of the maps
$$
{\mathcal O}_s(U,\mathcal{B}^{\rm r}(\mathcal{V})_\C)\ni F\mapsto F(T)\in\mathcal{B}^{\rm r}(\mathcal{V}),\,
{\mathcal O}_s(U)\ni f\mapsto f(T)\in\mathcal{B}^{\rm r}(\mathcal{V}),
$$
is clear. The second one is actually (unital and) multiplicative, which follows from the multiplicativiry of the usual analytic functional calculus  of $T$.

In fact, we have
$(Ff)(T)=F(T)f(T)$ for all   $F\in\mathcal{O}_s(U,  \B^{\bf r}(\mathcal{V})_\C),\, f\in{\mathcal O}_s(U)$. This follows from the equalities,
$$
(Ff)(T)=\frac{1}{2\pi i}\int_{\Gamma_0} F(\zeta)f(\zeta)(\zeta-T)^{-1}d\zeta=
$$
$$
\left(\frac{1}{2\pi i}\int_{\Gamma_0} F(\zeta)(\zeta-T)^{-1}d\zeta\right)
\left(\frac{1}{2\pi i}\int_{\Gamma} f(\eta)(\eta-T)^{-1}d\eta\right)=F(T)f(T),
$$
obtained as in the classical case (see \cite{DuSc}, Section VII.3), holding because $f$ is $\C$-valued and commutes with the operators in $ \B^{\bf r}(\mathcal{V})$. Here $\Gamma,\,\Gamma_0$ are the boundaries of two Cauchy
domains $\Delta,\,\Delta_0$ respectively, such that $\Delta\supset
\bar{\Delta}_0$, and $\Delta_0$ contains $\sigma_\C(T)$. 

Note that, in particular,  for every polynomial $P(\zeta)=\sum_{n=0}^m A_n\zeta^n$ with $A_n\in \B^{\bf r}(\mathcal{V})$  for all $n=0,1,\ldots,m$, we have  $P(T)=\sum_{n=0}^m A_n T^n\in  \B^{\bf r}(\mathcal{V})$ for all 
$T\in \B^{\bf r}(\mathcal{V})$.


\begin{Cor}\label{R_afc1}  Let $\V$ be a Banach $Cl$-space, let $U\subset\C$ be a conjugate symmetric open set, and let 
$T\in   \B^{\rm r}(\V)$, with $\sigma_\C(T)\subset U$.  There exists an assignment 
$$
{\mathcal O}_s(U,\Kk_n)\ni F\mapsto F(T)\in   \B^{\rm r}(\V),
$$ 
which is an $\R$-linear map, such that

(1) for all   $F\in\mathcal{O}_s(U,\Kk_n),\, f\in{\mathcal O}_s(U)$, we have $(Ff)(T)=F(T)f(T)$.

(2) for every polynomial $P(\zeta)=\sum_{n=0}^m \a_n\zeta^n,\,\zeta\in\C$, with $\a_n\in\Cc_n$  for all $n=0,1,\ldots,m$, we have  $P(T)=\sum_{n=0}^m \a_n T^n\in  \B^{\bf r}(\mathcal{V})$.  
\end{Cor}

{\it Proof.}  Note that the algebra ${\mathcal O}_s(U,\Kk_n)$, can  be regarded as a subalgebra of  the algebra
 ${\mathcal O}_s(U,   \B^{\rm r}(\V)_\C)$, whose elements are identified with left multiplication operators. Therefore, this corollary is a direct consequence of Theorem \ref{R_afc}.


\begin{Rem}\label{twofc}\rm The space $\mathcal{R}_s(\Omega,\Cc_n)$, introduced in Section 6,  can be independently defined, and it consists of the set of all $\Cc_n$-valued functions, which are {\it slice monogenic} in the 
sense of \cite{CoSaSt}, Definition 2.2.2 (or slice regular, as called in this work).  They are used in  \cite{CoSaSt} to define  a  functional calculus for tuples of not necessarily commuting real linear operators. Specifically, with a slightly modified notation, given an arbitrary family $(T_0,T_1,\ldots,T_n)$ , acting 
on the real space $\V$, it is associated with the operator ${\bf T}=\sum_{j=0}^n T_j\otimes e_j$, acting on the two-sided 
$\Cc_n$-module $\V_n=\V\otimes_\R\Cc_n$. In fact, the symbol ''$\otimes$'' may (and will) be omitted. Moreover, as alluded in \cite{CoSaSt}, page 83, we may work on a Banach $Cl$- space $\V$, and using operators from $\B^{\rm r}(\V)$.

 Roughly speaking, after fixing a Clifford  operator, each regular $\Cc_n$-valued function  defined in a neighborhood $\Omega$ of its $Cl$-spectrum is associated with another Clifford operator, replacing formally the paravector variable with that operator. This constraction is  explained in Chapter 3 of  \cite{CoSaSt}. 

 For an operator $T\in\mathcal{B}^{\rm r}(\mathcal{V})$, the {\it right $S$-resolvent} is defined via the formula

\begin{equation}\label{kqfc} 
S_R^{-1}(\s,T)=-(T-\s^*)(T^2-2\Re(\s)T+\Vert\s\Vert)^{-1},\,\,\s\in
\rho_{Cl}(T)
\end{equation}
(which is the right  version of formula (3.5) from \cite{CoSaSt}; see also formula (4.47) from \cite{CoSaSt}). Fixing an element $\kappa\in\S_n$, and a spectrally saturated open set $\Omega\subset\Pp_n$, for $\Phi\in\mathcal{R}_s(\Omega,\Cc_n)$ one sets 

\begin{equation}\label{qfc}
\Phi(T)=\frac{1}{2\pi}\int_{\Sigma_\kappa}\Phi(\s)d\s_\kappa S_R^{-1}(\s,T),
\end{equation}
 where $\Sigma_\kappa$  consists of a finite family of closed curves, piecewise smooth, positively oriented, being the boundary of the set  $\Theta_\kappa=\{\s=u+v\kappa\in\Theta;u,v\in\R\}$, where $\Theta\subset\Omega$ is a spectrally saturated open set containing  $\sigma_{Cl}(T)$, and $d\s_\kappa=-\kappa du\wedge dv$. Formula (\ref{qfc}) is a slight extension of the   (right)  functional calculus, as defined in \cite{CoSaSt}, Theorem 3.3.2 (see also formula (4.54) from \cite{CoSaSt}).

Our Corollary \ref{R_afc1} constructs, in particular, an analytic functional calculus with functions from ${\mathcal O}_s(U,\Kk_n)$, where $U$ is a  neighborhood of the complex spectrum  of a given Cliffordian operator, leading to another Clifford operator, replacing formally  the complex variable with that operator. We can show that those functional calculi are equivalent. This is a consequence of the 
isomorphism of the spaces  ${\mathcal O}_s(U,\Kk_n)$ and $\mathcal{R}_s(U_\sigma,\Cc_n)$, implied by Theorems \ref{C_afc} and \ref{equiv-ons-dom1}.

Let us give a direct argument concerning the equivalence of those analytic functional calculi. 
Because the space $\mathcal{V}_\C$ is also a $Cl$-space, we may apply these formulas to the extended
operator $T_\C\in\mathcal{B}^{\rm r}(\mathcal{V}_\C)$,  replacing $T$ by $T_\C$ in formulas (\ref{kqfc}) and (\ref{qfc}). In fact,
using the properties of the morphism $T\mapsto T_\C$ (see beginning of Section 7), we deduce that  $S_R^{-1}(\s,T)_\C=S_R^{-1}(\s,T_\C)$.

For the function  $\Phi\in\mathcal{R}_s(\Omega,\Cc_n)$ there exists a function $F\in{\mathcal O}_s(\Omega,\Kk_n)$ such that $F_\sigma=\Phi$,  by Theorem \ref{C_afc}. Denoting by
 $\Gamma_\kappa$ the boundary of a Cauchy domain in $\C$ containing the compact set
 $\cup\{\sigma(\s);\s\in\overline{\Theta_\kappa}\}$, we can write 
 $$
 \Phi(T_\C)=\frac{1}{2\pi}\int_{\Sigma_\kappa}\left(\frac{1}{2\pi i}\int_{\Gamma_\kappa}F(\zeta)(\zeta-\s)^{-1}d\zeta\right)
 d\s_\kappa S_R^{-1}(\s,T_\C)=
 $$
 $$
 \frac{1}{2\pi i}\int_{\Gamma_\kappa}F(\zeta)\left(\frac{1}{2\pi}\int_{\Sigma_\kappa}(\zeta-\s)^{-1}d\s_\kappa S_R^{-1}(\s,T_\C)\right)d\zeta.
 $$

It follows from the complex linearity of $S_R^{-1}(\s,T_\C)$, and via an argument  similar to that for getting 
formula (4.49) in \cite{CoSaSt}, that
$$
(\zeta-\s)S_R^{-1}(\s,T_\C)=S_R^{-1}(\s,T_\C)(\zeta-T_\C)-1,
$$
whence
$$
(\zeta-\s)^{-1}S_R^{-1}(\s,T_\C)=S_R^{-1}(\s,T_\C)(\zeta-T_\C)^{-1}+
(\zeta-\s)^{-1}(\zeta-T_\C)^{-1},
$$      
and therefore,
 
$$
\frac{1}{2\pi}\int_{\Sigma_\kappa}(\zeta-\s)^{-1}d\s_\kappa S_R^{-1}(\s,T_\C)=
 \frac{1}{2\pi}\int_{\Sigma_\kappa}d\s_\kappa S_R^{-1}(\s,T_\C)
(\zeta-T_\C)^{-1}+
$$
$$
\frac{1}{2\pi}\int_{\Sigma_\kappa}(\zeta-\s)^{-1}d\s_\kappa  
(\zeta-T_\C)^{-1}= (\zeta-T_\C)^{-1},
$$
because
$$
 \frac{1}{2\pi}\int_{\Sigma_\kappa}d\s_\kappa S_R^{-1}(\s,T_\C)=1\,\,\,{\rm and}
 \,\,\,\frac{1}{2\pi}\int_{\Sigma_\kappa}(\zeta-\s)^{-1}d\s_\kappa =0,
$$
as in Theorem 4.8.11 from \cite{CoSaSt}, since the $\Kk_n$-valued function $\s\mapsto(\zeta-\s)^{-1}$ is analytic in a neighborhood of the 
set $\overline{\Theta_\kappa}\subset \C_\kappa$ for each $\zeta\in\Gamma_\kappa$, respectively. Therefore $\Phi(T_\C)=\Phi(T)_\C=F(T_\C)=F(T)_\C$, implying $\Phi(T)=F(T)$.

Conversely,  choosing a function  $F\in{\mathcal O}_s(\Omega,\Kk_n)$,  and denoting by  $\Phi\in\mathcal{R}_s(\Omega,\Cc_n)$  its  Cauchy transform, the previous computation in reverse order shows that $\Phi(T)=F(T)$. Consequently, for a fixed 
$T\in\B^{\rm r}(\mathcal{V})$, the maps $\Theta:\mathcal{R}_s(\Omega,\Cc_n)\mapsto \B^{\rm r}(\mathcal{V})$, with $\Theta(\Phi)=
\Phi(T)$, and $\Psi:{\mathcal O}_s(\Omega,\Kk_n)\mapsto \B^{\rm r}(\mathcal{V})$, with $\Psi(F)=F(T)$, we must have the equality
$\Psi=\Theta\circ\C[*]$, where $C[*]$ is the Cauchy transform. 
\end{Rem}


\begin{Rem}\label{spmapth}\rm  Unlike in \cite{CoSaSt,CoGaKi}, our approach permits to obtain a version of the {\it spectral mapping theorem} in a classical stile, via direct arguments. Recalling that $\mathcal{R}_{s,n}(\Omega)$ is the subalgebra of  $\mathcal{R}_s(\Omega,\Cc_n)$ whose  elements are also in $\mathcal{IF}(\Omega,\Cc_n)$ (see Theorem \ref{C_afc}),
 for every operator  $T\in\mathcal{B}^{\rm r}(\mathcal{V})$ and every function $\Phi\in\mathcal{R}_{s,n}(\Omega)$  one has $\sigma_{Cl}(\Phi(T))=
\Phi(\sigma_{Cl}(T))$, via Theorem 3.5.9 from \cite{CoSaSt}. Using our approach, for every function $f\in\mathcal{O}_s(U)$, one has $f(\sigma_\C(T))=\sigma_\C(f(T))$, directlly from the corresponding (classical) spectral mapping theorem in \cite{DuSc}. This  result is parallel to that 
from \cite{CoSaSt} mentioned above, also giving an explanation for the former, via the isomorphism of the spaces 
$ \mathcal{O}_s(U)$ and $\mathcal{R}_{s,n}(\Omega)$
\end{Rem}


\section{Application to Tuples of Real Operators }

The special case exhibited at the beginning of  Remark \ref{twofc}  is largely treated in \cite{CoSaSt}, in connection with slice regular  functions and analytic functional calculus.
In this section, we shall briefly present some consequences of the results from the previous ones, valid, inparticular, for not necessarily commuting
tuples of linear operators acting on a given real Banach space $\V$.
We adapt, with our notation, the framework of \cite{CoSaSt}.
For a fixed integer $n\ge1$, we consider the real vector space
$\V_n=\V\otimes_\R\Cc_n$. The elements of $\V_n$ will be written
under the form ${\bf v}=\sum_J v_J e_J$, with $v_J\in\V$, where 
$J\prec\N_n$, and the symbol $"\otimes``$ will be omitted. The space
$\V_n$ is a  two-sided $\Cc_n$-module, with the operations
$$
(\sum_J u_Je_J)(\sum_K v_Ke_K)=\sum_{J,K}u_Jv_Ke_Je_K,
(\sum_K v_Ke_K)(\sum_J u_Je_J)=\sum_{K,J}u_Jv_Ke_Ke_J,
$$
for all elements $\sum_J u_Je_J\in\Cc_n,\sum_K v_Ke_K\in\V_n$.

Fixing a norm $\Vert*\Vert$ on $\V$, we define a norm on $\V_n$
by $\Vert{\bf v}\Vert_n=\sum_J\Vert v_J\Vert$, where ${\bf v}=
\sum_J v_Je_J$.

Following  \cite{CoSaSt}, the space  $\V_n$ is a Banach $\Cc_n$-module if there exists a constant $C\ge1$ such that 
$\Vert{\bf av}\Vert_n$, $\Vert{\bf va}\Vert_n$ are both majored by
$C\vert{\bf a}\vert\Vert{\bf v}\Vert_n$, for all  ${\bf a}\in\Cc_n$, ${\bf v}\in\V_n$. With our terminology, in this case the space
$\V_n$ is a Banach $Cl$-space. 

Let $\B(\V)$ be the algebra of $\R$-linear operators of the real 
Banach space $\V$. For a fixed family $\{T_J\}_{J\prec\N_n}$, we define an operator ${\bf T}=\sum_J T_Je_J$, acting on $\V_n$  via the formula
$$
{\bf T}({\bf v})=\sum_ J\sum_K  T_J(v_K )e_J e_K,\,\,{\bf v}=\sum_K v_K e_K\in\V_n.
$$
The set of all operators of this form will be denoted by $\B_n(\V_n)$.  Setting the norm $\Vert{\bf T}\Vert=\sum_J\Vert T_J\Vert$,
the set  $\B_n(\V_n)$ is a unital real Banach algebra (see \cite{CoSaSt}, page 82). 

Note that, with ${\bf T}\in\B_n(\V_n)$ represented as above, and ${\bf v}=\sum_K v_K e_K\in\V_n$, we have
$$
{\bf T}({\bf v\a})=\sum_{J,K,L}a_LT_J(v_K)e_Je_Ke_L={\bf T}({\bf v})\a,\,\, \a=\sum_L a_Le_L,
$$
showing that  ${\bf T}\in \B^{\rm r}(\V_n)$. In other words,  $\B_n(\V_n)\subset \B^{\rm r}(\V_n)$, and the inclusion is strict, as
simple examples show. An analytic functional calculus for the operator  ${\bf T}=\sum_J T_J e_J\in \B^{\rm r}(\V_n)$ can be obtained directly, as a consequence of Theorem  \ref{R_afc}. Nevertheless, as in  \cite{CoSaSt}, with minor modifications, we may replace the algebra 
 $\B^{\rm r}(\V_n)$ by the algebra   $\B_n(\V_n)$.

First of all, we consider the complexification $\V_{n\C}=\V_n+i\V_n$  of the real vector space $\V_n$. If 
${\bf T}=\sum_J T_Je_J\in\B_n(\V_n)$,
and ${\bf w}=\sum_K w_Ke_K\in \V_{n\C}$, with $w_K=u_K+iv_K, u_K,v_K\in\V_n$, then
$$
{\bf T}_\C({\bf w})=\sum_{J,K}(T_J(u_K)+iT_J(v_K))e_Je_K=\sum_{J,K}T_{J\C}(w_K)e_Je_K\in\B_n(\V_{n\C}).
$$

The conjugation ${\bf w}=\sum_K w_Ke_K\mapsto\bar{\bf w}=\sum_K\bar{w}_Ke_K$ on  $\V_{n\C}$, say $C$, induces a conjugation
${\bf S}\mapsto{\bf S}^\flat$ via the definition ${\bf S}^\flat=C{\bf S}C$ for all ${\bf S}\in\B_n(\V_{n\C})$. Moreover,
${\bf S}={\bf S}^\flat$ if and only if ${\bf S}(\V_n)\subset \V_n$, and ${\bf T}_\C^\flat={\bf T}_\C$ for all ${\bf T}\in\B_n(\V_n)$.

As in the case of real operators (see \cite{Vas3}), we define the complex spectrum of the operator ${\bf T}\in\B_n(\V_n)$  by 
the equality $\sigma_\C({\bf T})=\sigma(({\bf T}_\C)$.


\begin{Thm}\label{R_afc2} Let $\V$ be a real Banach space,  let ${\bf T}=\sum_J T_Je_J$ acting on $\V_n$, and 
 let $U\subset\C$ be a conjugate symmetric open set with $\sigma_\C({\bf T})\subset U$.  Then there exists an assignment 
$$
{\mathcal O}_s(U,  \B_n(\V_{n\C})\ni F\mapsto F({\bf T})\in  \B_n(\V_n)
$$ 
is an $\R$-linear map, and the restricted  map
$$
{\mathcal O}_s(U)\ni f\mapsto f({\bf T})\in  \B_n(\V_n)
$$
is a unital real algebra morphism.

Moreover, the following properties are true:

(1) for all   $F\in\mathcal{O}_s(U, \B_n(\V_{n\C}),\, f\in{\mathcal O}_s(U)$, we have $(Ff)({\bf T})=
F({\bf T})f({\bf T})$.

(2) for every polynomial $P(\zeta)=\sum_{n=0}^m A_n\zeta^n,\,\zeta\in\C$, with $A_n\in \B_n(\V_n)$  for all $n=0,1,\ldots,m$, we have  $P({\bf T})=\sum_{n=0}^m A_n {\bf T}^n\in \B_n(\V_n)$.  
\end{Thm}

{\it Proof.} The proof follows the lines of that of Theorem \ref{R_afc} (see also Theorem 2 from \cite{Vas3}).  We only note that a suitable version of the first part of Theorem \ref{afcro1} (or Theorem 1 and Remark 7 from \cite{Vas3})  is necessary, with minor modifications. The detaile are left to the reader.


\begin{Rem}\rm \label{tuples}  As in  \cite{CoSaSt}, fixing an arbitrary family $\{T_j\}_{j=1}^n$, one  considers the operator ${\bf T}=\sum_ {j=1}^n T_je_j$, which acts on $\V_n$  via the formula
$$
{\bf T}({\bf v})=\sum_ {j=1}^n\sum_K a_K T_j(v_K )e_j e_K,\,\,{\bf v}=\sum_K v_K e_K\in\V_n.
$$ 
In the previous statement we  may replace the algebra $\mathcal{O}_s(U, \B_n(\V_{n\C)}$ by the smaller algebra  $\mathcal{O}_s(U,\Kk_n)$, regarding its elements as left multiplication operators (see Corollary \ref{R_afc1}).  In this way, we obtain an analytic functional calculus for tuples of real operators, using analytic stem functions. 

A parallel version of some results from this section was  discussed in the first part of  \cite{CoSaSt}, using  the slice regular functions to  construct functional calculi.  
\end{Rem}

\end{document}